\documentclass[11pt]{article}
\usepackage[dvips]{graphics}
\input{amssym.def}
\input{amssym.tex}

\title{\Large\bf 
Two-dimensional  contact of two\\ different power-law graded elastic bodies}
\author{Y.A. Antipov$^a$, S.M.\ Mkhitaryan$^{b,c}$\\
$^a$Department of Mathematics, Louisiana State University,\\
Baton Rouge LA 70803, USA\\
$^b$Department of Mechanics of Elastic and Viscoelastic Bodies,\\
National Academy of Sciences,\\
$^c$Department of Mathematics and Physics,\\ 
 National University of Architecture and Construction\\
 Yerevan 0009, Armenia}

\newcommand{\R}{\mathop{\rm Re}\nolimits}
\newcommand{\const}{\mbox{const}}

\newcommand{\beqa}{\begin{eqnarray}}
\newcommand{\eeqa}[1]{\label{#1}\end{eqnarray}}
\newcommand{\bequ}{\begin{equation}}
\newcommand{\eequ}[1]{\label{#1}\end{equation}}
\newcommand{\Ga}{\alpha}
\newcommand{\Gb}{\beta}
\newcommand{\Gd}{\delta}

\newcommand{\Gve}{\varepsilon}
\newcommand{\Gf}{\phi}

\newcommand{\Gg}{\gamma}

\newcommand{\Gn}{\eta}

\newcommand{\Gt}{\theta}

\newcommand{\Gs}{\sigma}

\newcommand{\Gj}{\tau}

\newcommand{\Go}{\omega}
\newcommand{\Gx}{\xi}

\newcommand{\Gz}{\zeta}
\newcommand{\GD}{\Delta}
\newcommand{\GF}{\Phi}
\newcommand{\GG}{\Gamma}

\newcommand{\GS}{\Sigma}

\newcommand{\GY}{\Psi}

\newcommand{\CK}{{\cal K}}

\newcommand{\beq}{\begin{equation}}
\newcommand{\eeq}{\end{equation}}
\newcommand{\barr}{\begin{eqnarray}}
\newcommand{\earr}{\end{eqnarray}}
\newcommand{\beqn}{\begin{equation*}}
\newcommand{\eeqn}{\end{equation*}}
\newcommand{\barrn}{\begin{eqnarray*}}
\newcommand{\earrn}{\end{eqnarray*}}

\newcommand{\fr}{\frac}

\begin{document}
\maketitle

%\noindent
%{\it Abbreviated title:  Subsonic interfacial crack with a friction zone} 
%\renewcommand{\theequation}{\thesection.\arabic{equation}}

\begin{abstract}

Previous study of contact of power-law graded materials concerned the contact of a rigid body (punch) with
an elastic inhomogeneous foundation whose inhomogeneity is characterized by the Young modulus varying with depth as 
a power function. This paper models Hertzian and adhesive contact of two elastic inhomogeneous power-law graded 
bodies with different exponents. The problem is governed by an integral equation with two different power kernels.
A nonstandard method of Gegenbauer orthogonal polynomials for its solution is proposed. It leads to infinite system
of linear algebraic equations of a special structure. The integral representations of the system coefficients are evaluated, and the properties
of the system are studied. It is shown that if the exponents coincide, the infinite system admits a simple exact solution
that corresponds to the case when the Young moduli are different but the exponents are the same. Formulas for
the length of the contact zone, the pressure distribution, and the surface normal displacements of the contacting
bodies are obtain in the form convenient for computations. Effects of the mismatch in the Young moduli exponents
are studied. A comparative analysis of the Hertzian and adhesive contact models clarifies the effects of the
surface energy density on the contact pressure, the contact zone size, and the profile of the contacting bodies
outside the contact area.

\end{abstract}

Keywords: 

Two different power-law graded bodies

Hertzian contact

Adhesive contact

Novel method of Gegenbauer polynomials

\setcounter{equation}{0}

\section{Introduction} 

Interest in contact problems of interaction of  bodies with elastic inhomogeneous foundations
was originated in the forties of the previous century when civil engineers started taking into account
the inhomogeneity properties of soil foundations.  For the last thirty years, when
novel functionally graded materials (FGMs) were designed and the necessity of the study of their properties
arose (Saleh, 2020), this interest became even stronger. 

 One of the most interesting classes of FGMs comprises inhomogeneous materials 
whose modulus of elasticity $E$ varies with depth according to the power law, $E(z)=E_\Ga z^\Ga$.
The first approximate solution of a contact problem of an axisymmetric foundation with the modulus of elasticity
$E(z)=E_\Ga z^\Ga$  and subjected to a point force $P$ applied to the boundary was obtained by Klein (1955) in the form
\beq
\Gs_z=\fr{A P z^A}{2\pi R^{A+2}}, \quad 
\Gs_r=\fr{A P z^{A-2}r^2}{2\pi R^{A+2}},\quad \Gs_\Gt=0,
\quad \tau_{rz}=\fr{A P z^{A-1}r}{2\pi R^{A+2}},
 \label{1.1}
 \eeq
 where $R=\sqrt{r^2+z^2}$. This solution satisfies the equilibrium equations for any values of the constant $A$. However, in general, 
  the compatibility conditions for the strains are not met.
 It was found (Klein, 1955) that in only two cases, (1) $A=\Ga+3$, $\nu=(2+\Ga)^{-1}$ and  (2) $A=\Ga+2$, $\nu=(1+\Ga)^{-1}$,
 where $\nu$ is the Poisson ratio,
 the compatibility conditions are fulfilled. Also, in these particular cases it is possible to recover the normal displacement in the interior of the body
 by explicitly integrating the strain $e_z=(\Gs_z-\nu\Gs_r)(E_\Ga z^{\Ga})^{-1}$. Upon passing to the limit $z\to 0$ in the resulting
 formula for the displacement, this 
 gives  
  the normal displacement on the surface of the half-space,
 $w(x,y,0)=\Ga P/{\pi E_\Ga} r^{-\Ga-1}$, where $r=\sqrt{x^2+y^2}.$
Based on the solution obtained in these cases, Klein (1955) suggested  to extrapolate the formula for the displacement
$w(x,0)$, valid in only these two particular cases, to the general case when the Poisson ratio  $\nu$ and the exponent $\Ga$ are not connected by any
relation.  

 Leknitskii (1962)  considered the plane problem of a wedge with a variable modulus of elasticity.   On applying
the separation of variables method to the equilibrium equations he obtained an exact 
formula for the  stress $\Gs_r$ in a half-plane $\{|x|<\infty, y>0\}$ when $E=E_\Ga y^\Ga$ for any constant Poisson ratio $\nu$. 
By separating the variables in the equation for the Airy function Rostovtsev (1964) not only rederived the Lekhnitskii formula
for the stress but also obtained the exact representation for the normal displacement in the cases of concentrated and distributed
normal load applied to the boundary. In addition, he proved that in a general three-dimensional inhomogeneous
medium it is impossible to have a radial distribution of stresses. In particular, Rostovtsev (1964) showed that the Lekhnitskii problem,
when being axisymmetric and stated for a half-space with the modulus of elasticity $E(z)=E_\Ga z^\Ga$, 
except for the two particular cases examined by Klein (1955), does not have solutions with a radial distribution of stresses. 
  
In many contact problems, it is required to find  only the pressure distribution in the interior of the contact zone and the surface displacements in its exterior when the displacements in the contact area are prescribed. For such problems, when the experimental data
show that the classical elastic homogeneous isotropic half-space does not 
accurately model the deformable foundation, Korenev  (1960) introduced the concept of the kernel of a linearly deformable elastic foundation.  By means of the kernel of a foundation the normal
displacement may be expressed through the pressure distribution as
\beq
w(x,y)=\int\limits_{\Go} 
 K(x-\Gx,y-\Gn)p(\Gx,\Gn)d\Gx d\Gn, \quad (x,y)\in\Go,
 \label{1.6}
 \eeq
where $\Go$ is the contact area, and the kernel admits the representation 
$K(x-\Gx,y-\Gn)=\CK(r)$,  $r=\sqrt{(x-\Gx)^2+(y-\Gn)^2}$. 
A matrix generalization of Korenev's kernel of the foundation was proposed by Popov (1982) for the case when
the normal and tangential displacements in the contact area are prescribed, while the normal and tangential traction
components in the contact area are to be determined. 
The kernel of elastic homogeneous isotropic foundation is well-known, $\CK(r)=(1-\nu^2)(\pi E r)^{-1}$. 
Owing to the Klein's solution (1955) obtained for two particular cases, 
the kernel  $\CK(r)=(\pi D_\Ga)^{-1} r^{-\Ga-1}$ is often referred to as the kernel of
an  elastic inhomogeneous half-space whose Young modulus varies according to the power law, $E=E_\Ga z^\Ga$.
Korenev (1960) introduced five other kernels of linearly-deformable foundations. They are
$$
\CK_{I}(r)=\fr{A}{\sqrt{r^2+\Gd^2}},\quad \CK_{II}(r)=\fr{A}{2\Gd^2}\exp\left(-\fr{r^2}{4\Gd^2}\right), \quad 
\CK_{III}(r)=AK_0(\Gd r), 
$$
\beq
\CK_{IV}(r)=\fr{A}{2\pi r}e^{-\Gd r}, \quad \CK_V(r)=\fr{A}{2\pi(r^2+\Gd^2)},
\label{1.7}
\eeq
where $A$ and $\Gd$ are positive parameters determined by tests, and $K_0(\Gd r)$ is the modified Bessel function. 
Note that  the axisymmetric contact problem of a circular stamp indented into an elastic half-space characterized by the kernel 
$\CK_{IV}(r)$ was solved in terms of spheroidal functions by Mkhitaryan (2015). Recently, Antipov and Mkhitaryan (2021)
analyzed bending of a strip-shaped and a half-plane-shaped plate 
lying on an elastic foundation
characterized by the kernel $\CK_{III}(r)$. 
 
The majority of work on plane and axisymmetric contact problems of power-law graded materials  concern  the indentation of
a rigid two-dimensional or axisymmetric stamp  into a half-plane or a half-space. In the case of a single contact zone, the plane problem reduces
to the integral equation
\beq
\Gg_0\int_{-b}^b\fr{p(\Gx)d\Gx}{|x-\Gx|^\Ga}=\Gd-f(x), \quad -b<x<b,\quad 0<\Ga<1,
\label{1.9}
\eeq
where $\Gd$ is the indentation of the stamp,  the function $f(x)$ describes the stamp profile, and $\Gg_0$ is a function of $\Ga$.
The solution of this equation in the class of functions admitting integrable singularities at the endpoints $\pm b$
exists and unique. It can be constructed by a variety of methods including the method of Abelian integrals (see for example,
Gakhov, 1966),
the method of dual integral equations, the Wiener-Hopf method, and the method of orthogonal polynomials. The solution
of this integral equation by the last method is presented in Section 5 of this paper.
Popov (1967) considered the more advanced case of this plane problem when there are two separate contact zones.
He  reduced the problem to two separately solvable equations with the Weber-Schafheitlin kernel and solved them
approximately by the method of the Jacobi polynomials.
The first exact solutions to the axisymmetric case were obtained  by the method of dual integral equations (Korenev, 1957; Mossakovskii, 1958) under the assumption
of the frictionless contact of a stamp and a power-law graded foundation.  The same problem was later solved (Popov, 1961)
by the Wiener-Hopf method. The method of Abelian operators was  applied  by Popov (1973) to derive 
an exact solution to the axisymmetric
problem of non-slipping adhesive contact of a punch  with a power-law graded elastic half-space.

During the last twenty five years plane and axisymmetric contact problems  of a stamp and a half-plane and a half-space
with the Young modulus $E=E_\Ga z^\Ga$ have become the subject of interest (Giannakopoulos and Suresh, 1997a, 1997b;
Giannakopoulos and Pallot, 2000; Chen et al, 2009a, 2009b, Guo, 2011; Willert, 2018; Jin et al, 2021)
 due to modeling of micro- and nano-indentation processes arising in nanotechnology and therefore the necessity of characterization
of mechanical properties of a variety of biological materials with
sizes approaching molecular or atomic dimension (Guo et al, 2011). These authors considered the Johnson-Kendal-Roberts 
(JKR) adhesive 
 model (Johnson et al, 1971; Johnson, 1985) to examine plane and axisymmetric contact of a rigid punch with a half-plane
and half-space, respectively, when the Young modulus of the foundation varies with depth according to a power-law.
The feature of the JKR model is that it admits integrable singularities of the contact pressure at the endpoints and determines
the contact zone length (radius) from the condition of minimum of the total energy. The total energy $U_{total}$ is defined
to be a sum
of the elastic strain energy $U_e$ and the loss of surface energy $U_s$.  Another approach to modeling of adhesive
contact, the  Maugis-Dugdale model (Maugis, 1992) was recently employed (Jin et al, 2021)
to examine axisymmetric contact of a punch and a power-law graded half-space. This model assumes that 
the cohesive stress is constant within
 the cohesive zone outside the contact area.

There have been relatively limited efforts in studying Hertzian and adhesive contact 
of two elastic bodies whose Young moduli are power-functions of depth.
Popov and Savchuk (1971) considered the axisymmetric Hertzian model of contact of two bodies having different Young moduli  
$E_1(z)=e_1 z^\Ga$  and $E_2(z)=e_2 (-z)^\Ga$ but the same exponents. They also took 
 into account the surface effects according to the Shtayerman (1949) model. Power-law kernels arise
 in the problem of  computing equilibrium measures
for problems with attractive-repulsive kernels of the form $K(x-y)=\Ga^{-1}|x-y|^\Ga-\Gb^{-1}|x-y|^\Gb$ Cutleb et al, (2021). 
 For this problem, they
 proposed a numerical method of recursively generated banded and approximately banded operators acting 
on expansions in ultraspherical polynomial bases.
 To the best of the authors knowledge, neither two-dimensional nor axisymmetric problem of Hertzian or JKR adhesive contact of two elastic bodies with different Young moduli,
$E_1(z)=e_1 z^{\Ga_1}$ and $E_2(z)=e_2 (-z)^{\Ga_2}$, have been considered in the literature.

In this paper we aim to
analyze the plane contact problem of two different power-law graded bodies.
In Section 2, we formulate the problem and reduce it  to the integral equation with two kernels of the form
\beq
\int_{-b}^b\left(\fr{\Gg_1}{|x-\Gx|^{\Ga_1}}+\fr{\Gg_2}{|x-\Gx|^{\Ga_2}}\right)p(\Gx)d\Gx=\Gd-f(x), \quad -b<x<b,
\label{1.10}
\eeq
where $\Gd$ is a rigid body displacement to be determined from an equilibrium condition, $p(x)$ is the 
pressure distribution, $\Gg_1$ and $\Gg_2$ are some positive parameters,  $f(x)=f_1(x)+f_2(x)$, $y=f_1(x)$ 
and $y=-f_2(x)$ are the profiles of 
the contacting bodies,  $0<\Ga_2<\Ga_1<1$. This equation may be interpreted as a full integral equation with a single power kernel $|x-\Gx|^{-\Ga_1}$ with the second kernel serving as a regular part
(Gakhov, 1966). However, the method of Abelian operators, when applied, leads to a Fredholm integral 
equation whose kernel is a chain of singular integrals, and does not produce the solution in the form convenient for numerical
purposes. 

In Section 3, we describe the method of solution that 
expands  the unknown function $ p(bt)$ in terms
of the Gegenbauer polynomials $C_n^{\Ga_1/2}(t)$ with weight $(1-t^2)^{(\Ga_1-1)/2}(t)$ and reduces the task
of finding the expansion coefficients to solution of an infinite system of linear algebraic coefficients with coefficients represented by integrals
possessing the polynomials  $C_n^{\Ga_1/2}(t)$ and $C_m^{\Ga_2/2}(t)$.  We manage to evaluate
these integrals. The coefficients have certain remarkable properties which substantially simplify the system.
We also show that in the limit case $\Ga_2\to\Ga_1$, the solution of the infinite system can be derived explicitly, and it coincides with the solution 
of the contact problem of two bodies with different  power-law Young moduli and the same exponent, $E_1=e_1y^{\Ga}$
and $E_2=e_2 (-y)^{\Ga}$.

In Section 4, we derive formulas for the length of the contact zone,  the parameter $\Gd$, the pressure distribution, and the normal displacement on the surface outside the contact zone in the 
form  convenient for computations. 
We emphasize that all the formulas except for the displacement are free of integrals.  
We also discuss the results of numerical tests.

In Section 5, we derive a closed-form solution of the  problem of Hertzian contact of two bodies whose moduli of elasticity have the same exponents
$\Ga_1=\Ga_2=\Ga$
but different factors $e_1$ and $e_2$. We obtain exact formulas not only for the contact zone length and the pressure but also
for the normal displacement outside the contact area.

In Section 6, we analyze the JKR model for both cases, when $\Ga_1=\Ga_2$ and $\Ga_1>\Ga_2$.  In both cases we compute
the elastic strain energy and the total energy. In the former 
case we obtain a transcendental equation for the contact zone half-length $b$ and show that it is possible  to pass to the limit
as $\Ga_j\to 0$.  In the case $\Ga_1>\Ga_2$ we derive the equation for $b$ approximately by computing the derivative
of the strain energy numerically.  We show that in both cases the solution to the JKR model coincides with the solution to 
the Hertzian model when the surface energy half-density $\Gg_s\to 0$.

\setcounter{equation}{0}

\section{Formulation}\label{form}

\begin{figure}[t]
\centerline{
\scalebox{0.6}{\includegraphics{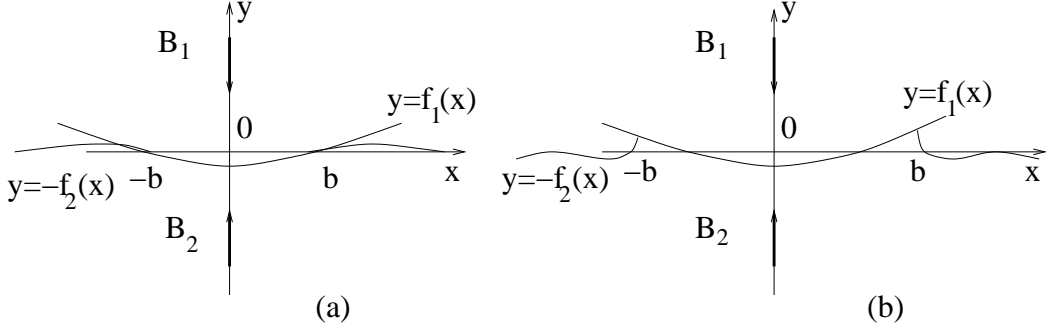}}
}
\caption{Contact of two different power-law graded elastic bodies with Young moduli $E_1(y)=e_1y^{\Ga_1}$ (body $B_1$) 
and $E_2(y)=e_2(-y)^{\Ga_2}$ (body $B_2$). (a): Hertzian model and (b): adhesive JKR model.}
\label{fig1}
\end{figure}

The problem of interest is the one of modeling of two-dimensional 
contact of two inhomogeneous elastic bodies, $B_1$ and $B_2$ 
(Figure 1 (a)). The lower surface of the upper body  $B_1$ and the upper surface of the lower body $B_2$ are
described by curves $y=f_1(x)$ and $y=-f_2(x)$. The functions
$f_1(x)$ and $f_2(x)$ are even, continuously differentiable and share the tangent line
$y=0$ at the point $x=0, y=0$,  the origin of the Cartesian coordinates $(x,y)$, that is $f_1(0)=f_2(0)=0$ and 
$f'_1(0)=f'_2(0)=0$. The bodies are inhomogeneous
whose Poisson ratios $\nu_1$ and $\nu_2$ are constant, while the Young moduli vary according
to a power law and equal $E_1(y)=e_1y^{\Ga_1}$ and $E_2(y)=e_2(-y)^{\Ga_2}$, respectively, where $e_1$ and $e_2$
are positive constants, $0<\Ga_j<1$, $j=1,2$.
The bodies are subjected to compression by forces applied to the bodies parallel to the $y$-axis
with the resultant force $P$  balanced by the contact pressure $p(x)$ arising in the contact area $(-b,b)$,
and the parameters  $b$ is unknown {\it a priori}. 
 We also assume that the curve $y=f_1(x)$ is  convex upward, while the second curve
$y=f_2(x)$ is either convex downward or flat or at least locally convex upward.
To proceed with the contact modeling, we make the following Hertzian assumptions:

\begin{itemize}

\item the contact area is significantly less than the bodies sizes,

\item the friction is absent, and the only nonzero traction component is $\Gs_y=-p(x)$, where $p(x)$ is the normal pressure,

\item the normal and tangential elastic displacements in the contact area are significantly smaller than the contact zone length..

\end{itemize}

Following Shtayerman (1949) we write the vertical displacements of any two points $A_1\in B_1$ and $A_2\in B_2$
which, as a result of compression,  become the same point, a point $A$. These displacements
are $f_1(x-u_1)+v_1-\Gd_1$ and   $-f_2(x+u_2)-v_2+\Gd_2$. Here, $(u_1,v_1)$ and $(-u_2,-v_2)$ are the elastic displacements
of the points $A_1$   and $A_2$, and the constants $\Gd_1$ and $\Gd_2$ are forward displacements of distant points.  Approximating 
$f_1(x-u_1)\approx f_1(x)$ and $f_2(x+u_2)\approx f_2(x)$, we can write  at the point of contact $A$
\beq
v_1+v_2=\Gd-f_1(x)-f_2(x), \quad -b<x< b, \quad \Gd=\Gd_1+\Gd_2.
\label{2.1}
\eeq
 The parameter $\Gd$
is to be determined {\it a posteriori} from the condition
\beq
\int_{-b}^b p(x)dx=P.
\label{2.2}
\eeq
 
We next use the Rostovtsev relation (Rostovtsev, 1964, p.747) between the normal displacement and contact pressure for a half-plane
to write down the displacements $v_1$ and $v_2$ in the contact area
  \beq
  v_j(x)=\fr{\Gt_j}{\Ga_j}\int_{-b}^b\fr{p(\Gx)d\Gx}{|x-\Gx|^{\Ga_j}}, \quad -b<x<b, \quad j=1,2.
  \label{2.3}
  \eeq
 Here,
$$
\Gt_j= \fr{C_j(1-\nu_j^2)q_j}{(\Ga_j+1)e_j}\sin\fr{\pi q_j}{2}, 
\quad q_j=\sqrt{(1+\Ga_j)\left(1-\fr{\Ga_j\nu_j}{1-\nu_j}\right)},
$$
\beq
C_j=\fr{2^{\Ga_j+1}}{\pi\GG(\Ga_j+2)}
\GG\left(\fr{\Ga_j}{2}-\fr{q_j}{2}+\fr32\right)\GG\left(\fr{\Ga_j}{2}+\fr{q_j}{2}+\fr32\right).
\label{2.4}
\eeq
Substituting the integral representations of the displacements $v_j$ into the condition $(\ref{2.1})$ we derive
the governing integral equation for the contact pressure distribution $p(x)$
\beq
\int_{-b}^b\left(\fr{\Gt_1}{\Ga_1|x-\Gx|^{\Ga_1}}+\fr{\Gt_2}{\Ga_2|x-\Gx|^{\Ga_2}}\right)p(\Gx)d\Gx=\Gd-f_1(x)-f_2(x), \quad -b<x<b.
\label{2.5}
\eeq
To show that this equation gives rise to the integral equation of Hertzian contact of two homogeneous elastic bodies,
we first rewrite the equation in the form
\beq
\int_{-b}^b\left[\fr{\Gt_1}{\Ga_1}(|x-\Gx|^{-\Ga_1}-1)+\fr{\Gt_2}{\Ga_2}(|x-\Gx|^{-\Ga_2}-1)\right]p(\Gx)d\Gx=\Gd_0-f_1(x)-f_2(x),
 \quad -b<x<b,
\label{2.5'}
\eeq
where
\beq
\Gd_0=\Gd-\left(\fr{\Gt_1}{\Ga_1}+\fr{\Gt_2}{\Ga_2}\right)P
\label{2.6}
\eeq
is a free constant. Then, 
by letting $\Ga_j\to 0$, $j=1,2$,  and taking into account that
\beq
\lim_{\Ga_j\to 0}\fr{|x-\Gx|^{-\Ga_j}-1}{\Ga_j}=\ln\fr{1}{|x-\Gx|},  
\label{2.7}
\eeq
and also that
\beq
q_j\to 1, \quad C_j\to \fr{2}{\pi}, \quad \Gt_j\to\Gt_j^\circ=\fr{2(1-\nu_j^2)}{\pi E_j}\quad {\rm as} \;\Ga_j\to 0,
\label{2.8}
\eeq
we obtain the classical integral equation when $E_j(y)=E_j=\const$ (Shtayerman, 1949)
\beq
(\Gt_1^\circ+\Gt_2^\circ)\int_{-b}^b \ln\fr{1}{|x-\Gx|} p(\Gx)d\Gx=\Gd_0-f_1(x)-f_2(x), \quad -b<x<b.
\label{2.9}
\eeq

\setcounter{equation}{0}
  
\section{Solution of the integral equation}\label{int} 

To solve the integral equation (\ref{2.5}), it will be convenient to rewrite it in the interval $(-1,1)$
\beq
\int_{-1}^1\left(\fr{A_1}{|t-\tau|^{\Ga_1}}+\fr{A_2}{|t-\tau|^{\Ga_2}}\right)p(b\tau)d\tau=\Gd-f(bt), \quad -1<t<1,
\label{3.3}
\eeq
where
\beq
A_j=\fr{\Gt_jb^{1-\Ga_j}}{\Ga_j}.
\label{3.4}
\eeq
The right-hand side of equation (\ref{3.3}) possesses the unknown parameter $\Gd$. To eliminate it from the equation, we represent
the function $p(bt)$ as 
\beq
p(bt)=\Gf^{(1)}(t)+\Gd\Gf^{(2)}(t)
\label{3.5}
\eeq
and deduce
\beq
\int_{-1}^1\left(\fr{A_1}{|t-\tau|^{\Ga_1}}+\fr{A_2}{|t-\tau|^{\Ga_2}}\right)\Gf^{(j)}(\tau)d\tau=g^{(j)}(t), \quad -1<t<1, \quad j=1,2.
\label{3.6}
\eeq
where $g^{(1)}(t)=-f(bt)$, $g^{(2)}(t)=1$. The equilibrium condition (\ref{2.2}) expresses the unknown parameter $\Gd$
through the solutions $\Gf_1$ and $\Gf_2$ of the equations (\ref{3.6})
which share the kernel and have different right-hand sides. We have
\beq
\Gd=\left(\fr{P}{b}-\int_{-1}^1\Gf^{(1)}(\Gj)d\Gj\right)\left(\int_{-1}^1\Gf^{(2)}(\tau)d\tau\right)^{-1}.
\label{3.7}
\eeq

\subsection{Infinite system of algebraic equations}\label{inf}  

Without loss of generality we assume further that $\Ga_1>\Ga_2$ and denote
\beq
\Gb_{n}(\Ga)=\fr{\pi(\Ga)_n}{n!\cos\fr{\pi\Ga}{2}}, \quad n=0,1,\ldots,
\label{3.8}
\eeq
where $(\Ga)_n=\Ga(\Ga+1)\ldots(\Ga+n-1)$ is the factorial symbol.
Owing to the spectral relation for the Gegenbauer polynomials
\beq
\int_{-1}^1 \fr{C_n^{\Ga/2}(\Gj)d\Gj}{|t-\Gj|^{\Ga}(1-\Gj^2)^{(1-\Ga)/2}}=\Gb_{n}(\Ga)C_n^{\Ga/2}(t), \quad -1<t<1,
\quad 0<\Ga<1,
\label{3.9}
\eeq
and the orthogonality property of these polynomials
\beq
\int_{-1}^1 C_n^{\Ga/2}(t) C_m^{\Ga/2}(t)(1-t^2)^{(\Ga-1)/2}dt=h_n(\Ga) \Gd_{mn}, \quad m,n=0,1,\ldots,
\label{3.10}
\eeq
we seek the solution in the series form 
\beq
\Gf^{(j)}(t)=(1-t^2)^{(\Ga_1-1)/2}\sum_{n=0}^\infty\GF_{n}^{(j)} C_n^{\Ga_1/2}(t), \quad -1<t<1, \quad j=1,2.
\label{3.11}
\eeq
Here, $\GF_n^{(j)}$ are unknown coefficients,  $\Gd_{mn}$ is the Kronecker symbol, and
\beq
h_n(\Ga)=\fr{\pi 2^{1-\Ga}\GG(n+\Ga)}{n!(n+\fr{\Ga}{2})\GG^2(\fr{\Ga}{2})}.
\label{3.12}
\eeq

In the  integral equations of a rigid stamp indented into an inhomogeneous power-law graded
half-plane or Hertzian contact of two bodies with $\Ga_1=\Ga_2$, there is only one power-law
kernel. In these particular  cases, the series coefficients can be derived explicitly by substituting  the expansion (\ref{3.11})
into the integral equation and taking into account the spectral relation (\ref{3.9}) and the orthogonality
property (\ref{3.10}).  In contrast to this, when $\Ga_1\ne \Ga_2$, 
 we have the second term in the kernel,  and, in general, the series coefficients 
cannot be found exactly.
On substituting (\ref{3.11}) into (\ref{3.6}) we have
\beq
A_1\sum_{n=0}^\infty\Gb_n(\Ga_1)\GF_{n}^{(j)}C_n^{\Ga_1/2}(t)+
A_2\sum_{n=0}^\infty\GF_{n}^{(j)}
\int_{-1}^1
\fr{G_n(\Gj)(1-t^2)^{(\Ga_2-1)/2}d\Gj}
{|t-\Gj|^{\Ga_2}}
=g^{(j)}(t), \quad -1<t<1,
\label{3.13}
\eeq
where $G_n(\Gj)=C_n^{\Ga_1/2}(\Gj)(1-t^2)^{(\Ga_1-\Ga_2)/2}$. Since $\Ga_1>\Ga_2$, we may expand the
functions $G_n(\Gj)$ in terms of the Genebauer polynomials $C_m^{\Ga_2/2}(\Gj)$
\beq
G_n(\Gj)=\sum_{m=0}^\infty G_m^{(n)} C_m^{\Ga_2/2}(\Gj), \quad -1<\Gj<1.
\label{3.14}
\eeq
According to the orthogonality relation (\ref{3.10}) the coefficients of the expansion are found to be
\beq
G_m^{(n)}=\fr{H_m^{(n)}}{h_m(\Ga_2)}, \quad H_m^{(n)}=\int_{-1}^1C_n^{\Ga_1/2}(\Gj)C_m^{\Ga_2/2}(\Gj)(1-\Gj^2)^{(\Ga_1-1)/2}d\Gj.
\label{3.15}
\eeq
Notice that $H_m^{(n)}=0$ if $m<n$. Indeed, the degree-$m$ polynomial $C_m^{\Ga_2/2}(\Gj)$ is a linear combination
of the monomials $1,\Gj, \ldots, \Gj^m$ or, equivalently, a linear combination of the Gegenbauer polynomials 
$C_0^{\Ga_1/2)}(\Gj), C_1^{\Ga_1/2}(\Gj),\ldots, C_m^{\Ga_1/2)}(\Gj)$, and by the orthogonality relation (\ref{3.10})
$H_m^{(n)}=0$  provided $m<n$. Now, if we substitute the series (\ref{3.14}) back to equation (\ref{3.13}),
use the spectral relation (\ref{3.10}) for the Gegenbauer polynomials $C_m^{\Ga_2/2}(\Gj)$ and change the order of summation,
we find
\beq
A_1\sum_{n=0}^\infty\Gb_n(\Ga_1)\GF_n^{(j)}C_n^{\Ga_1/2}(t)+
A_2\sum_{n=0}^\infty\GY_{n}^{(j)}\Gb_n(\Ga_2)C_n^{\Ga_2/2}(t)
=g^{(j)}(t), \quad -1<t<1,
\label{3.16}
\eeq
where
\beq
\Psi_{n}^{(j)}=\sum_{m=0}^n G_n^{(m)}\GF_{m}^{(j)}.
\label{3.17}
\eeq
The equation (\ref{3.16}) can be recast by using the orthogonality relation (\ref{3.10}) and written as an infinite system
of algebraic equations. We have
\beq
A_1\Gb_n(\Ga_1)h_n(\Ga_1)\GF_{n}^{(j)}+A_2\sum_{m=0}^\infty\GY_{m}^{(j)}\Gb_m(\Ga_2)H_m^{(n)}=g^{(j)}_n, \quad n=0,1,\ldots,
\label{3.18}
\eeq
where
\beq
g_n^{(j)}=\int_{-1}^1 g_j(t)C_n^{\Ga_1/2}(t)(1-t^2)^{(\Ga_1-1)/2}dt.
\label{3.19}
\eeq
It is possible to simplify the system derived.
On changing the order of summation in the series in the system (\ref{3.18}) and using the relations (\ref{3.15})
and (\ref{3.17}) we obtain
\beq
\sum_{m=0}^\infty \GY_{m}^{(j)}\Gb_m(\Ga_2)H_m^{(n)}=\sum_{m=0}^\infty L_{nm}\GF_{m}^{(j)},
\label{3.20}
\eeq
where
\beq
L_{nm}=\sum_{k=\max(m,n)}^\infty \fr{H_k^{(n)}H_k^{(m)}\Gb_k(\Ga_2)}{h_k(\Ga_2)}, \quad m,n=0,1,\ldots,
\label{3.21}
\eeq
and therefore the system has the form
\beq
A_1\Gb_n(\Ga_1)h_n(\Ga_1)\GF_{n}^{(j)}+A_2\sum_{m=0}^\infty L_{nm}\GF_{m}^{(j)}
=g_{n}^{(j)}, \quad n=0,1,\ldots.
\label{3.22}
\eeq

\subsection{Evaluation of the integrals $H_k^{(m)}$}\label{H}

We remind that $H_k^{(m)}=0$ if $k<m$. To compute the integrals (\ref{3.15}) when $k\ge m$, we use the formula
$$
\int_{-1}^1(1-x)^\Ga(1+x)^{\nu-1/2}C_m^\mu(x)C_n^\nu(x)dx=\fr{2^{\Ga+\nu+\fr12}\GG(\Ga+1)\GG(\nu+\fr12)\GG(\nu-\Ga+n-\fr12)}
{m! n! \GG(\nu-\Ga-\fr12)\GG(\nu+\Ga+n+\fr32)}
$$
\beq
\times\fr{\GG(m+2\mu)\GG(n+2\nu)}{\GG(2\mu)\GG(2\nu)}
{}_4 F_3\left(\begin{array}{cccc}
-m, & m+2\mu,& \Ga+1,&\Ga-\nu+\fr32\\
\mu+\fr12, &\nu+\Ga+n+\fr32,&\Ga-\nu-n+\fr32;&1\\
\end{array}
\right),
\label{3.24}
\eeq
where  $\R\Ga>-1$,  $\R\nu>-\fr12$ and ${}_4 F_3$ is the generalized hypergeometric function. This relation can be derived from the general formula for the Jacobi polynomials
 (Bateman and Erdelyi, 1954, formula 16.4(20)).
Notice that the corresponding formulas for the Gegenbauer polynomials (Bateman and Erdelyi, 1954, formula 16.3(16))
and Gradshteyn and Ryzhik, 1994, formula 7.314(7)) have the same error: instead of $\GG(\nu+\Ga+n+\fr32)$ 
in the right-hand side in (\ref{3.24}) they write $\GG(\nu-\Ga+n+\fr32)$.

 On adjusting the relation (\ref{3.24})
to our case when $k\ge m$
we have 
\beq
H_k^{(m)}=\fr{2\sqrt{\pi}(-1)^m \GG(\fr{\Ga_1+1}{2})(\Ga_2)_k}{m!\GG(\fr{\Ga_1}{2})(m+\Ga_1)}\GS,
\label{3.25}
\eeq
where
\beq
\GS=\sum_{l=m}^k
\fr{(-1)^l(\Ga_2+k)_l (\fr{\Ga_1+1}{2})_l}{(k-l)!(l-m)!(\Ga_1+m+1)_l(\fr{\Ga_2+1}{2})_l}.
\label{3.25'}
\eeq
This sum can be evaluated and the formula for $H_k^{(m)}$ simplified. We 
make the substitution $l-m=i$,
use the property of the factorial symbols
\beq
(a)_{m+i}=(a+m)_i(a)_m, \quad (k-n)!=\fr{(-1)^n k!}{(-k)_n},\quad k\ge n,
\label{3.26}
\eeq
and express the sum $\GS$ through the function  ${}_3 F_2$
\beq
\GS=\fr{(-1)^m(\fr{\Ga_1+1}{2})_m(\Ga_2+k)_m}{(k-m)!(\Ga_1+m+1)_m(\fr{\Ga_2+1}{2})_m}
{}_3 F_2\left(\begin{array}{ccc}
-k+m, & \Ga_2+k+m, & \fr{\Ga_1+1}{2}+m\\
\fr{\Ga_2+1}{2}+m, & \Ga_1+2m+1; & 1\\
\end{array}
\right).
\label{3.27}
\eeq
For the generalized hypergeometric function ${}_3 F_2$ in the right-hand side we can employ Whipple's formula (Wipple, 1925)
\beq
{}_3 F_2\left(\begin{array}{ccc}
a, & b, & c\\
\fr{a+b+1}{2}, & 2c; & 1
\end{array}
\right)=\fr{\sqrt{\pi}\GG(c+\fr12)\GG(\fr{a+b+1}{2})\GG(\fr{1-a-b}{2}+c)}
{\GG(\fr{a+1}{2})\GG(\fr{b+1}{2})\GG(\fr{1-a}{2}+c)\GG(\fr{1-b}{2}+c)}
\label{3.28}
\eeq
and obtain the following representation for $\GS$:
$$
\GS=\fr{(-1)^m(\fr{\Ga_1+1}{2})_m(\Ga_2+k)_m}{(k-m)!(\Ga_1+m+1)_m(\fr{\Ga_2+1}{2})_m}
$$
\beq
\times\fr{\sqrt{\pi}\GG(\fr{\Ga_1}{2}+m+1)\GG(\fr{\Ga_2+1}{2}+m)\GG(\fr{\Ga_1-\Ga_2}{2}+1)}
{\GG(\fr{m-k+1}{2})\GG(\fr{\Ga_1+m+k}{2}+1)\GG(\fr{\Ga_2+m+k+1}{2})\GG(\fr{\Ga_1-\Ga_2+m-k}{2}+1)}.
\label{3.29}
\eeq
This formula implies that $\GS=0$  and therefore $H_k^{(m)}=0$ if $k=m+1+2l$, $l=0,1,\ldots$.
In the case when $k-m$ is even, $k=m+2l$, $l=0,1,\ldots$, we substitute (\ref{3.29}) into (\ref{3.25}) and find
$$
H_{m+2l}^{(m)}=
\fr{2\sin\fr{\pi(\Ga_2-\Ga_1)}{2}
\GG(\fr{\Ga_1+1}{2})(\fr{\Ga_1+1}{2})_m
(\fr{\Ga_1}{2})_{m+1}
\GG(\fr{\Ga_1-\Ga_2}{2}+1)\GG(\fr{\Ga_2+1}{2})}
{\pi m!\GG(\Ga_2)(m+\Ga_1)(\Ga_1+m+1)_m}
$$
\beq
\times
\fr{\GG(2l+2m+\Ga_2)\GG(l+\fr12)\GG(\fr{\Ga_2-\Ga_1}{2}+l)}{(2l)!\GG(\fr{\Ga_1}{2}+m+l+1)\GG(\fr{\Ga_2+1}{2}+m+l)}.
\label{3.30}
\eeq
 On exploiting further the properties of the $\GG$-function it is possible to give to formula (\ref{3.30}) 
a different form
\beq
H_{m+2l}^{(m)}=\fr{\sqrt{\pi}\GG(\fr{\Ga_1+1}{2})}{\GG(\fr{\Ga_1}{2}+1)}
\fr{(\Ga_1)_m(\Ga_2/2)_m}{m!(\fr{\Ga_1}{2}+1)_m}
\fr{(\fr{\Ga_2-\Ga_1}{2})_l(\fr{\Ga_2}{2}+m)_l}{l!(\fr{\Ga_1}{2}+m+1)_l}.
\label{3.31}
\eeq
This formula is simpler and convenient for analysis of the coefficients asymptotics as $l\to\infty$. 
Taking into account the asymptotic relation
\beq
\fr{\GG(z+a)}{\GG(z+b)}\sim z^{a-b}, \quad z\to \infty, \quad |\arg z|<\fr{\pi}{2},
\label{3.32}
\eeq
we derive
\beq
H_{m+2l}^{(m)}\sim C_m l^{\Ga_2-\Ga_1-2}, \quad l\to \infty,
\label{3.33}
\eeq
where $C_m$ are constants.

Having computed the coefficients $H_k^{(m)}$ we consider now two cases, $n=0,1,\ldots,m-1$
and $n=m,m+1,\ldots$ and evaluate the coefficients $H_k^{(n)}$.
In the former case according to formula (\ref{3.21}) and since $H_k^{(m)}=0$ if $k=m+2l+1$, $l=0,1,\ldots$,
we need to evaluate $H_k^{(n)}$ for $k=m+2l$ only. On replacing $m$ by $n$ and $k$ by $m+2l$ in (\ref{3.25})
and (\ref{3.29}) we should have
$H_{m+2l}^{(n)}=0$,  if  $n-m$  is  odd and   $l=0,1,\ldots.$
Otherwise, if $n-m$ is even,
\beq
H_{m+2l}^{(n)}=\fr{2^{\Ga_2-\Ga_1}\sqrt{\pi}\GG(\fr{\Ga_2+1}{2})\GG(\Ga_1+n)}{\GG(\fr{\Ga_1}{2})\GG(\Ga_2)\GG(\fr{\Ga_2-\Ga_1}{2})n!}
\fr{\GG(\fr{\Ga_2+m+n}{2}+l)\GG(\fr{m-n+\Ga_2-\Ga_1}{2}+l)}{\GG(\fr{m-n}{2}+l+1)\GG(\fr{\Ga_1+m+n}{2}+l+1)}, 
 \quad  l=0,1,\ldots,
 \label{3.34}
 \eeq
and their asymptotics  for large $l$ is the same as for $H_{m+2l}^{(m)}$. We have 
\beq
H_{m+2l}^{(n)}\sim C'_{mn} l^{\Ga_2-\Ga_1-2}, \quad l\to \infty,
\label{3.34'}
\eeq
where $C_{mn}$ are constants. 
We also give  another, more convenient for numerical purposes, representation of the coefficients $H_{m+2l}^{(n)}$ 
when $n-m$ is even
\beq
H_{m+2l}^{(n)}=\fr{\sqrt{\pi}\GG(\fr{\Ga_1+1}{2})}{\GG(\fr{\Ga_1}{2}+1)}
\fr{(\Ga_1)_n}{n!}
\fr{(\fr{\Ga_2}{2})_{(m+n)/2}}
{(\fr{\Ga_1}{2}+1)_{(m+n)/2}}
\fr{(\fr{\Ga_2-\Ga_1}{2})_{(m-n)/2}}{(\fr{m-n}{2})!}
\fr{(\fr{\Ga_2+m+n}{2})_l (\fr{m-n+\Ga_2-\Ga_1}{2})_l}
{(\fr{m-n}{2}+1)_l(\fr{\Ga_1+m+n}{2}+1)_l}.
 \label{3.35}
 \eeq

\subsection{Solution of the infinite system}\label{sol} 

By introducing new notations we rewrite the system (\ref{3.22}) in the canonical form
\beq
\GF_{n}^{(j)}+\Gg\sum_{m=0}^\infty R_{nm}\GF_{m}^{(j)}=d_{n}^{(j)}, \quad n=0,1,\ldots, \quad j=1,2,
\label{3.36}
\eeq
where
\beq
\Gg=\fr{A_2}{A_1}, \quad R_{nm}=\fr{L_{nm}}{\Gb_n(\Ga_1)h_n(\Ga_1)}, \quad d_{n}^{(j)}=\fr{g_{n}^{(j)}}{A_1\Gb_n(\Ga_1)h_n(\Ga_1)}.
\label{3.37}
\eeq
Owing to the fact that $H_{m+2l}^{(n)}=0$,  if  $n-m$  is  odd and   $l=0,1,\ldots.$, from formula (\ref{3.21})
we deduce that $L_{nm}=0$ and therefore $R_{nm}=0$ if $n-m$ is odd. We have also derived that 
$H_{k}^{(m)}=0$ if $k=m+2l+1$ and $l=0,1,\ldots$. This brings us to the following formulas for the coefficients
$L_{nm}$ when $m-n$ is even:
$$
L_{nm}=\sum_{l=0}^\infty H_{m+2l}^{(m)}H_{m+2l}^{(n)}\GD_{m+2l}, \quad n=0,1,\ldots,m-1,
$$
\beq
L_{nm}=\sum_{l=0}^\infty H_{n+2l}^{(n)}H_{n+2l}^{(m)}\GD_{n+2l}, \quad n=m,m+1,\ldots,
\label{3.38}
\eeq
where 
\beq
\GD_k=\fr{\GG(\fr{\Ga_2}{2})\GG(\fr{1-\Ga_2}{2})(k+\fr{\Ga_2}{2})}{\sqrt{\pi}},
\label{3.39}
\eeq
$H_{n+2l}^{(m)}$ is obtained by interchanging  $n$ and $m$ in  (\ref{3.34}), while  $H_{n+2l}^{(n)}$ will
coincide with (\ref{3.31}) if $m$ is replaced by $n$. To sum up, for all  $n,m=0,1,\ldots$, $L_{nm}=L_{mn}\ne 0$ if $n-m$ is even and $L_{nm}=0$  otherwise.

 Remark that owing to the asymptotic relations  (\ref{3.33}) and (\ref{3.34'}) and formula (\ref{3.39}) the coefficients
 in the series (\ref{3.38}) behave for large $l$ as $l^{2(\Ga_2-\Ga_1)-3}$ $(\Ga_1>\Ga_2)$, and therefore the series representations (\ref{3.38})
 for the coefficients $L_{nm}$ rapidly converge.
 
 On passing to the limit $\Ga_2\to \Ga_1$ we can show that the matrix of the infinite system is diagonal, the system
admits an exact solution that coincides with that associated with the contact problem of two bodies with the same 
exponent $\Ga_1=\Ga_2$.  Indeed, when $\Ga_1=\Ga_2$ from (\ref{3.31}) and (\ref{3.35}) we deduce that  in either case, $l>0$ or $n\ne m$, 
the coefficients $H_{n+2l}^{(m)}$ and $H_{m+2l}^{(n)}$ are equal to zero, and  the only nonzero 
coefficients are $H_n^{(n)}$. They are given by
\beq
H_n^{(n)}=\fr{\sqrt{\pi}\GG(\fr{\Ga_1+1}{2})(\Ga_1)_n}{\GG(\fr{\Ga_1}{2})(\fr{\Ga_1}{2}+n)n!}.
\label{3.40}
\eeq
This gives a simple formula for the coefficients $L_{nm}$. It is $L_{nm}=[H_n^{(n)}]^2\GD_n\Gd_{nm}$, and from (\ref{3.37}),
$R_{nm} =\Gd_{nm}$. The system (\ref{3.36}) has a diagonal matrix, and 
the coefficients  $\GF_{n}^{(j)}=(1+\Gg)^{-1}d_n^{(j)}$ are the same as those obtained by solving the integral equation (\ref{3.6})
when $\Ga_1=\Ga_2$ on using the standard method of orthogonal polynomials.

In the general case, when $0<\Ga_2<\Ga_1<1$, the infinite system (\ref{3.36}) does not admit an exact
solution. Its approximate solution is found by the reduction method. The off-diagonal elements of the 
matrix of the system $\Gd_{mn}+\Gg R_{mn}$ rapidly decay, and the numerical method demonstrates a rapid convergence.

The right-hand sides of the system (\ref{3.36}) are represented by the integrals (\ref{3.19}). The integral $g_n^{(2)}$ is evaluated immediately,
$g_n^{(2)}=\GG_0\Gd_{n0}$, where
\beq
\GG_0=\fr{\sqrt{\pi}\GG(\fr{\Ga_1+1}{2})}{\GG(\fr{\Ga_1}{2}+1)}.
\label{3.41}
\eeq
The other integral $g_n^{(1)}$ can be
computed  explicitly if we know the coefficients $a_k$ of the expansion
of the function $ f(bt)$ in terms of the Gegenbauer polynomials
\beq
 f(bt)=\sum_{k=0}^\infty a_k C_k^{\Ga_1/2}(t).
\label{3.42}
\eeq
These coefficients are always computed exactly if the function $f(bt)$ is a polynomial. 
Otherwise, we can employ either its approximate polynomial representation or use the corresponding  Gauss' quadrature formula. 
In the polynomial case, when all the coefficients $a_k=0$, $k>N$,  we apply the orthogonality property (\ref{3.10})
to  find $g_n^{(1)}=-a_n h_n(\Ga_1)$, $n=0,1,\ldots,N$, and $g_n^{(1)}=0$, $n>N$.

\setcounter{equation}{0}
  
\section{Solution of the contact problem}\label{cont}

\subsection{Parameter $\Gd$,  the contact zone $(-b,b)$, the contact pressure $p(x)$, and the normal displacements $v_j$}

After the system (\ref{3.36}) for the two right-hand sides $d_n^{(1)}$ and $d_n^{(2)}$ has been solved and the 
values of the  coefficients $\GF_n^{(1)}$ and  $\GF_n^{(1)}$ have been found we write down 
the series representations (\ref{3.11}) of the solutions $\Gf^{(1)}(t)$ and $\Gf^{(2)}(t)$ of the integral equations (\ref{3.6}). On substituting these series into (\ref{3.7})
we can express the unknown parameter $\Gd$ through the coefficients $\GF_0^{(1)}$ and  $\GF_0^{(2)}$
\beq
\Gd=\fr{P/b- \GF_0^{(1)}\GG_0}{\GF_0^{(2)}\GG_0}.
\label{4.5}
\eeq
On having this parameter we can write down the contact pressure as
\beq
p(x)=\Gf^{(1)}\left(\fr{x}{b}\right)+\Gd \Gf^{(2)}\left(\fr{x}{b}\right).
\label{4.6}
\eeq
Notice that the parameter $\Gg=\Ga_1\Gt_2(\Ga_2\Gt_1)^{-1}b^{\Ga_1-\Ga_2}$ and the right-hand sides of the system (\ref{3.36}) depend on the unknown parameter $b$. That is why the contact pressure also depends on this parameter.
Because of the smoothness of the bodies profiles the contact pressure has to be bounded at the points $x=\pm b$, $y=0$. 
Owing to the representations (\ref{3.11}) this  implies that the contact pressure vanishes at these points,
\beq
\lim_{t\to 1}[\Gf^{(1)}(t)+\Gd\Gf^{(2)}(t)]=0.
\label{4.7}
\eeq
Equivalently, this reads
\beq
\sum_{n=0}^\infty \fr{(\Ga_1)_n}{n!}\left(\GF_n^{(1)}+\fr{P/b-\GF_0^{(1)}\GG_0}{\GF_0^{(2)}\GG_0}\GF_n^{(2)}\right)=0.
\label{4.8}
\eeq
This is a transcendental equation with respect to the parameter $b$. On having solved this equation we can determine 
the parameter $\Gd$ and the contact pressure by formulas (\ref{4.5}) and (\ref{4.6}), respectively.

The final quantities we wish to determine are the displacements  $v_j(x)$ of the surface points outside the contact zone.
We assume that  the curvatures of the surfaces of interest are sufficiently small. Since  formula (\ref{2.3}) for the
normal displacement is valid not only in the contact area but also outside, we can write
 \beq
  v_j(tx)=A_j\int_{-1}^1\fr{p(b\tau)d\tau}{|\tau-t|^{\Ga_j}}, \quad |t|>1, \quad j=1,2.
  \label{4.9}
  \eeq
Using formula (\ref{3.5}) and substituting the series representations (\ref{3.11}) into (\ref{4.9}) we 
 write the displacements as follows:
\beq
v_j(x)=A_j\sum_{n=0}^\infty
[\GF_n^{(1)}+\Gd\GF_n^{(2)}] I_{n}\left(\fr{x}{b};\Ga_j\right), \quad |t|>1,
\label{4.10}
\eeq
where
\beq
I_{n}(t;\Ga_j)=\int_{-1}^1\fr{(1-\Gj^2)^{(\Ga_1-1)/2} C_n^{\Ga_1/2}(\Gj)d\Gj}{|\Gj-t|^{\Ga_j}}.
\label{4.11}
\eeq
Series representations of this integral are derived in Appendix A. Since the function $f(x)$ is even,
all coefficients $\GF_{2m+1}^{(j)}=0$, $m=0,1,\ldots,$ $j=1,2$, and therefore 
\beq
v_j(x)=A_j\sum_{n=0}^\infty
[\GF_{2n}^{(1)}+\Gd\GF_{2n}^{(2)}] I_{2n}\left(\fr{x}{b};\Ga_j\right), \quad |t|>1.
\label{4.11'}
\eeq
On differentiating these functions we find out that the derivatives $v_j'(x)$ are bounded at the points $x=\pm b$
if and only if the condition (\ref{4.8}) is satisfied. In other words, if the contact zone parameter $b$
is fixed by solving the transcendental equation (\ref{4.8}), then not only the pressure $p(x)$ vanishes at the endpoints
but also the profiles of the contacting bodies are smooth at the endpoints.

\subsection{Numerical results}

\begin{figure}[t]
\centerline{
\scalebox{0.6}{\includegraphics{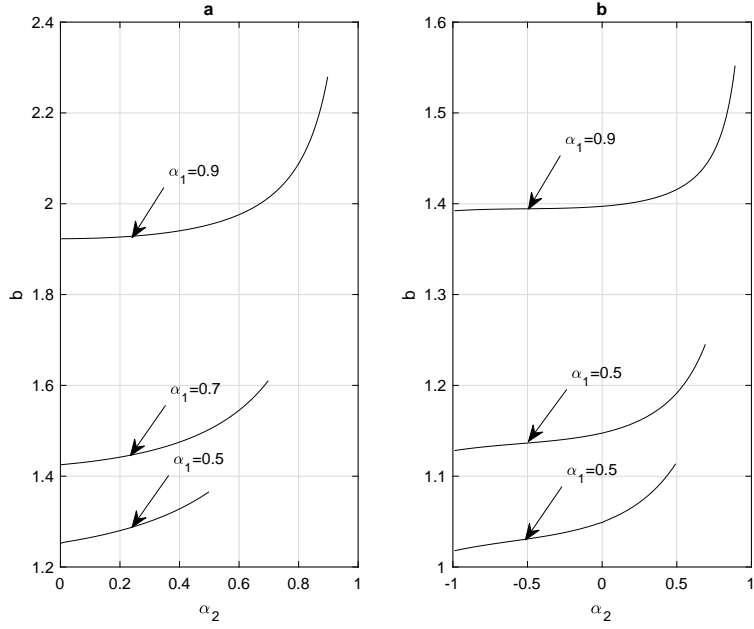}}}
\caption{
The  half-length $b$ of the contact zone $(-b,b)$  versus the parameter $\Ga_2\in(0,\Ga_1)$ 
for $\Ga_1=0.5$, $\Ga_1=0.7$ and $\Ga_1=0.9$ when (a) $f(x)=x^2$  ($Q_0=1$, $Q_1=0$) and
(b) $f(x)=x^4$  ($Q_0=0$, $Q_1=1$).}
\label{fig2}
\end{figure}

\begin{figure}[t]
\centerline{
\scalebox{0.6}{\includegraphics{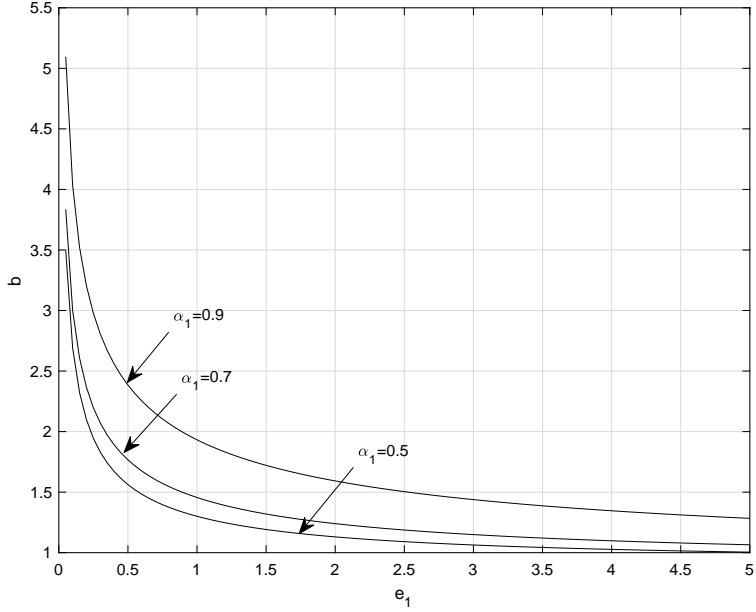}}}
\caption{
The  half-length $b$ of the contact zone $(-b,b)$  versus the parameter $e_1\in(0,5)$ 
for $\Ga_1=0.5$, $\Ga_1=0.7$ and $\Ga_1=0.9$ when $e_2=1$, $\Ga_2=0.3$, $f(x)=x^2$.}
\label{fig3}
\end{figure}

\begin{figure}[t]
\centerline{
\scalebox{0.6}{\includegraphics{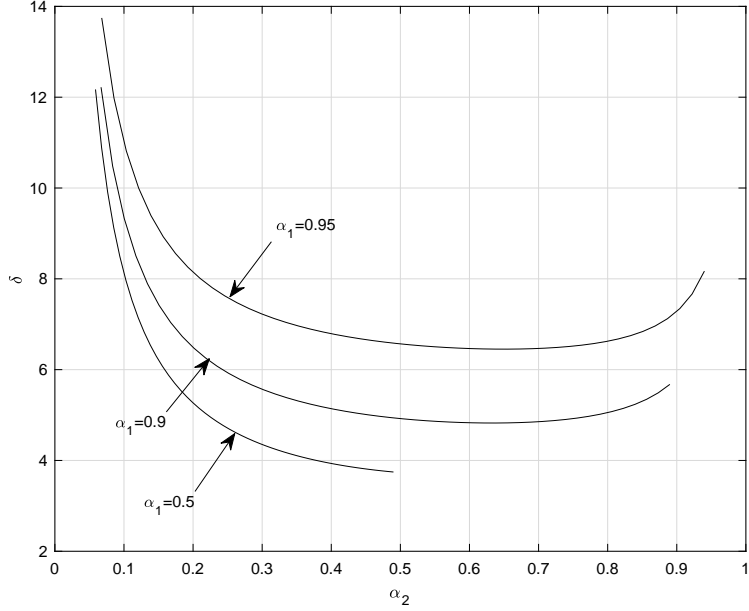}}}
\caption{
The  parameter $\Gd$ versus the exponent $\Ga_2\in(0,\Ga_1)$ 
for $\Ga_1=0.5$, $\Ga_1=0.9$ and $\Ga_1=0.95$ when  $f(x)=x^2$.}
\label{fig4}
\end{figure}

\begin{figure}[t]
\centerline{
\scalebox{0.6}{\includegraphics{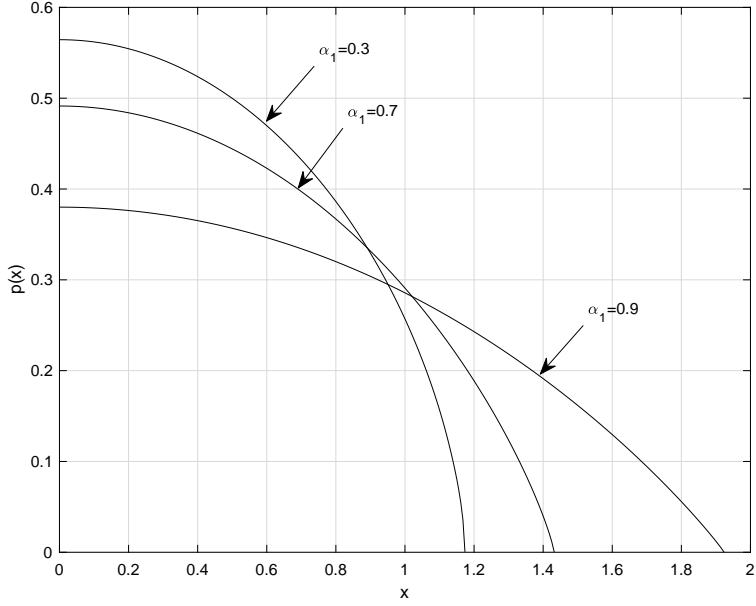}}}
\caption{
The contact pressure $p(x)$,  $x\in[0,b]$  for  $\Ga_2=0.1$ when
$\Ga_1=0.3$, $\Ga_1=0.7$, and $\Ga_1=0.9$ in the case  $f(x)=x^2$.}
\label{fig5}
\end{figure} 

\begin{figure}[t]
\centerline{
\scalebox{0.6}{\includegraphics{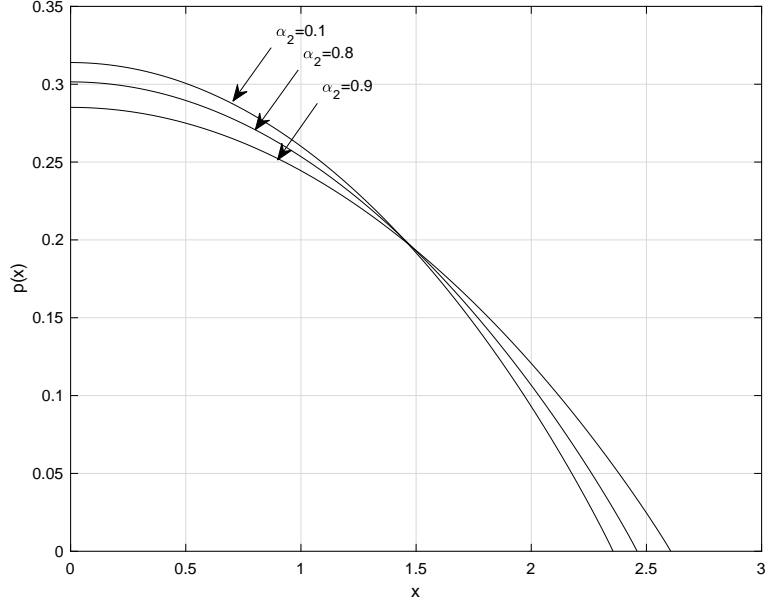}}}
\caption{
The contact pressure $p(x)$,  $x\in[0,b]$  for  $\Ga_1=0.95$ when
$\Ga_2=0.9$, $\Ga_2=0.8$, and $\Ga_2=0.1$ in the case  $f(x)=x^2$.}
\label{fig6}
\end{figure}

\begin{figure}[t]
\centerline{
\scalebox{0.6}{\includegraphics{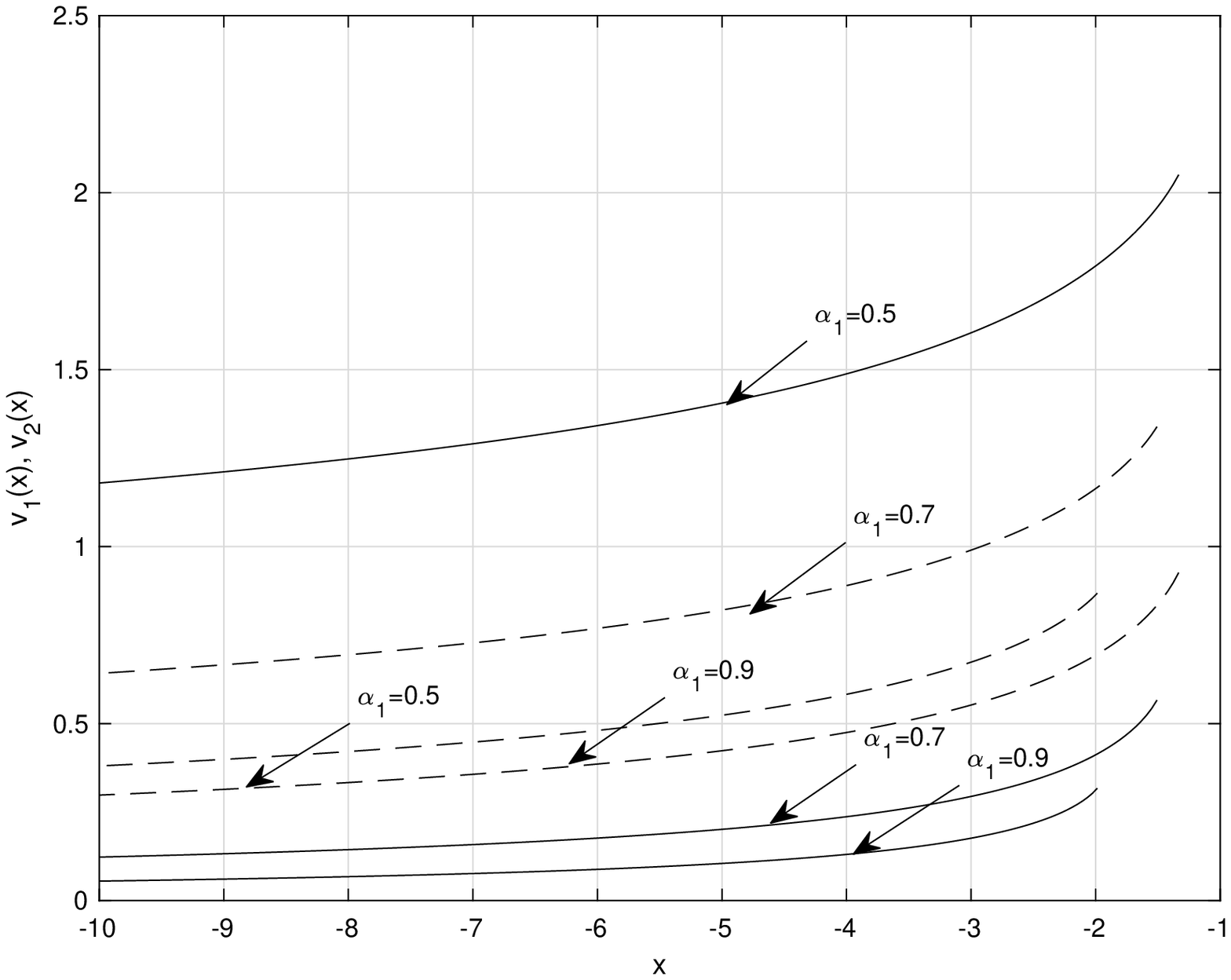}}}
\caption{
The normal displacements  $v_1(x)$ (solid curves) and $v_2(x)$ (broken curves),  $x\in[-10,-b]$  for  $\Ga_2=\Ga_1/2$ when
$\Ga_1=0.5$, $\Ga_1=0.7$, and $\Ga_1=0.9$ in the case  $f(x)=x^2$.}
\label{fig7}
\end{figure}

The functions $f_1(x)$  and $f_2(x)$ have to be continuously differentiable  and satisfy the conditions
$f_j(0)=f'_j(0)=0$, $j=0,1$. In the symmetric case, when both of the functions are even, in a neighborhood of the point $x=0$,
\beq
f_j(x)=\fr{f_j''(0)}{2}x^2+\fr{f^{(IV)}(0)}{24}x^4+\ldots, \quad j=1,2.
\label{4.12}
\eeq
For numerical tests, we confine ourselves to two polynomial cases of the function $f(x)=f_1(x)+f_2(x)$. They are

(1) $f(x)=Q_0x^2$, $Q_0>0$, and 

(2) $f(x)=Q_0x^2+Q_1x^4$. 

Case (1) occurs when one of the bodies has a parabolic profile,  while the second one is either flat or
also has a parabolic profile. In case (2), the profiles of the bodies are described by the polynomials  $f_j(x)=c_{0j} x^2+c_{1j}x^4$
with some real coefficients $c_{0j}$ and $c_{1j}$ chosen such that 
$Q_0=c_{01}+c_{02}\ge 0$ and $Q_1=c_{11}+c_{12}> 0$.

In case (1), we express the function $f(x)$ through the degree-0 and 2 Gegenbauer polynomials and have
\beq
f(bt)=\fr{b^2Q_0}{\Ga_1(\fr{\Ga_1}{2}+1)}\left[\fr{\Ga_1}{2}C_0^{\Ga_1/2}(t)+C_2^{\Ga_1/2}(t)\right].
\label{4.13}
\eeq
The orthogonality property (\ref{3.10}) yields $g_n^{(1)}=0$ for all $n$ unless $n=0$ or $n=2$. In these cases
\beq
g_0^{(1)}=-\fr{b^2Q_0\GG_0}{\Ga_1+2}, \quad  g_2^{(1)}=-\fr{b^2Q_0\sqrt{\pi}\GG(\fr{\Ga_1+3}{2})}{(\fr{\Ga_1}{2}+1)(\fr{\Ga_1}{2}+2)\GG(\fr{\Ga_1}{2})},
\label{4.14}
\eeq
where $\GG_0$ is given by (\ref{3.41}).

In case (2), the corresponding representation of the function $f(bt)$ has the form
$$
f(bt)=\fr{24Q_1 b^4 C_4^{\Ga_1/2}(t)}{\Ga_1(\Ga+2)(\Ga_1+4)(\Ga_1+6)}
$$
\beq
+\fr{2b^2C_2^{\Ga_1/2}(t)}{\Ga_1(\Ga_1+2)}\left(Q_0+\fr{6Q_1b^2}{\Ga_1+6}\right)+\fr{b^2C_0^{\Ga_1/2}(t)}{\Ga_1+2}
\left(Q_0+\fr{3Q_1b^2}{\Ga_1+4}\right),
\label{4.15}
\eeq 
Except for $g_0^{(1)}$, $g_2^{(1)}$, and $g_4^{(1)}$, all the terms $g_n^{(1)}$ equal 0. The nonzero terms are given by
$$
 g_0^{(1)}=-\fr{b^2\GG_0}{\Ga_1+2}\left(Q_0+\fr{3Q_1b^2}{\Ga_1+4}\right),
 $$
 $$
g_2^{(1)}=  -\fr{b^2\sqrt{\pi}\Ga_1\GG(\fr{\Ga_1+3}{2})}{2\GG(\fr{\Ga_1}{2}+3)}\left(Q_0+\fr{6Q_1b^2}{\Ga_1+6}\right),
$$
\beq
g_4^{(1)}=  -\fr{b^4\sqrt{\pi}\Ga_1(\Ga_1+2)\GG(\fr{\Ga_1+5}{2})}{4\GG(\fr{\Ga_1}{2}+5)}Q_1.
\label{4.16}
\eeq

For the numerical tests to be discussed we choose the resultant force and the Poisson ratios to be   $P=1$, $\nu_1=\nu_2=0.3$,
the resultant moment to be zero and the function $f(x)=f_1(x)+f_2(x)$  to be even. This choice gives rise to a solution 
symmetric with respect to the $y$-axis. Figure 2 presents the half-length $b$ of the contact zone for different values of the 
exponents $\Ga_1$ and $\Ga_2$ in the Young moduli of the bodies, $E_1=e_1y^{\Ga_1}$ and $E_2=e_2(-y)^{\Ga_2}$
when $e_1=e_2=1$ and (a) $f(x)=x^2$ and (b) $f(x)=x^4$ . It is seen that when $\Ga_1$ is fixed and $\Ga_2$ increases in the interval
$(0,\Ga_1)$, the contact zone length is also increasing. The same is true in the case  when $\Ga_2$ is fixed and $\Ga_1$ grows.
 On comparing the results presented in Figures 2 (a) and 2 (b) we see that  when the curvatures
of the contacting bodies profiles is decreasing the contact zone length is increasing.

Curves in Figure 3 give a clear demonstration of the dependence of the length of the contact zone
upon one of the factors $e_1$ and $e_2$ while the second one is kept fixed. The parameters for this diagram are chosen as
$e_2=1$,  $f(x)=x^2$,  $\Ga_2=0.3$, and $\Ga_1$ is equal to either $0.5$, $0.7$, or $0.9$. 

Figure 4 shows how the parameter $\Gd$ depends on  the second body exponent $\Ga_2\in(0,\Ga_1)$
when the exponent $\Ga_1$ is fixed and chosen to have the values $0.5$, $0.9$, and $0.95$.
The other parameters are $e_1=e_2=1$, and the function $f(x)=x^2$. It is seen that the parameter $\Gd$ increases 
as $\Ga_2\to 0$ and also when when $\Ga_1\to 1$ and $\Ga_2\to \Ga_1$.

The results of calculations of the pressure distribution $p(x)$ are shown in Figures 5 and 6. In both cases,  $e_1=e_2=1$, and  $f(x)=x^2$.
In Figure 5, the smaller exponent $\Ga_2$ is fixed as $\Ga_2=0.1$, while $\Ga_1$ is equal to either $0.3$, $0.7$, or $0.9$.
The corresponding values of the half-length $b$ of the contact zone are  computed to be $1.17365$,
$1.43214$, and  $ 1.92390$. The contact pressure $p(x)$ vanishes at the endpoints $\pm b$ of the contact zone and attains
its maximum at the origin. As the parameter $\Ga_1$ is increasing, the pressure maximum is decreasing. 
When the bigger exponent $\Ga_1$ is fixed (in Figure 6, $\Ga_1=0.95$), while the smaller exponent varies in the interval
$(0,\Ga_1)$, the variation of the pressure distribution  $p(x)$  for a fixed $x$ is not large (Figure 6).

The normal elastic displacements $u_y(x,0)=v_1(x)$ and $u_y(x,0)=-v_2(x)$ of the upper and lower elastic bodies
outside the contact zone are shown in Figure 8 (the displacements of the lower body $B_2$ are demonstrated by
broken curves). As before,  $e_1=e_2=1$, and  $f(x)=x^2$ and the functions 
$v_1(x)$ and $v_2(x)$ are even. For computations, we choose $\Ga_1$ to be either $0.5$, $0.7$, or $0.9$,
while $\Ga_2=\Ga_1/2$. Both displacements attain their maximum at the points $\pm b$. It has been
numerically verified that 
\beq
\lim_{x\to\pm b^\pm} [v_1(x)+v_2(x)]=\Gd-f_1(b)-f_2(b)
\label{4.17}
\eeq
that is consistent with the boundary condition (\ref{2.1}). The corresponding
values of the half-length of the contact zone are $b= 1.28951$ for $\Ga_1=0.5$,
$b=  1.46450$ when $\Ga_1=0.7$, and $b=1.94635$  in the case $\Ga_1=0.9$.

\setcounter{equation}{0}
  
\section{Hertzian contact of two power-law graded bodies when $E_1=e_1y^{\Ga}$
and $E_2=e_2 (-y)^{\Ga}$ }

Assume that the contacting bodies  $B_1$ and $B_2$ have the Young moduli $E_1=e_1 y^\Ga$ and $E_1=e_2 (-y)^\Ga$.
The governing equation (\ref{2.5}) with two kernels reduces to 
\beq
A\int_{-1}^1\fr{p(b\tau)d\tau}{|\tau-t|^{\Ga}}=\Gd-f_1(bt)-f_2(bt), \quad -1<t<1,
\label{5.1}
\eeq
where $A=\Ga^{-1}(\Gt_1+\Gt_2)b^{1-\Ga}$, and $\Gt_j$ are defined by (\ref{2.4}) with $\Ga_1=\Ga_2=\Ga$. The pressure distribution has to
be an even function, and we represent 
 the solution in the form
\beq
p(bt)=(1-t^2)^{(\Ga-1)/2}\sum_{n=0}^\infty\GF_{2n}C_{2n}^{\Ga/2}(t), \quad -1<t<1.
\label{5.2}
\eeq
On substituting this function into (\ref{5.1}), using the spectral relation (\ref{3.9}) and orthogonality property (\ref{3.10}) we find the coefficients
$\GF_{2n}$
\beq
\GF_{2n}=\fr{g_{2n}^{(1)}+\Gd\GG_0\Gd_{n0}}{A\Gb_{2n}(\Ga)h_{2n}(\Ga)},
\label{5.3}
\eeq
where
\beq
g_{2n}^{(1)}=-\int_{1}^1 C_{2n}^{\Ga/2}(t)(1-t^2)^{(\Ga-1)/2}f(bt)dt
\label{5.4}
\eeq
and $f(x)=f_1(x)+f_2(x)$. To determine the parameter $\Gd$, we  satisfy the equilibrium condition (\ref{2.2}) and obtain
\beq
\Gd=\fr{1}{\GG_0}\left(
-g_0^{(1)}+\fr{\pi AP}{b\cos\fr{\pi\Ga}{2}}\right).
\label{5.5}
\eeq
The half-length $b$ of the contact zone is the positive root of the following transcendental equation (our numerical tests reveal
that such a root is unique):
\beq
\sum_{n=0}^\infty \fr{(\Ga)_{2n}}{(2n)!}\GF_{2n}=0
\label{5.6}
\eeq
that reduces to
\beq
\fr{\Gd\GG_0\Ga}{2\GG(\Ga)}+\sum_{n=0}^\infty \fr{g_{2n}^{(1)}(2n)!(2n+\fr{\Ga}{2})}{\GG(\Ga+2n)}=0.
\label{5.7}
\eeq
The normal surface displacements of the bodies $B_1$ and $B_2$ are expressed through the integral $I_{2n}(x/b;\Ga)$
\beq
v_j(x)=A_j\sum_{n=0}^\infty \GF_{2n} I_{2n}\left(\fr{x}{b};\Ga\right), \quad |x|>b,\quad j=1,2,
\label{5.8}
\eeq
where
\beq
A_j=\fr{\Gt_jb^{1-\Ga}}{\Ga}, \quad I_{2n}(t;\Ga)=\int_{-1}^1\fr{(1-\Gj^2)^{(\Ga-1)/2} C_{2n}^{\Ga/2}(\Gj)d\Gj}{|\Gj-t|^{\Ga}}.
\label{5.9}
\eeq
This integral is a particular case $\Ga_j=\Ga$ of the integral $I_n(t;\Ga_j)$ evaluated in Appendix A
and given by (A.5) and (A.6).
As in the case $\Ga_1>\Ga_2$, the displacements $v_j$ and their first derivative are bounded as $x\to \pm b^\pm$,
and the contacting surfaces are smooth at the endpoints.
For numerical purposes, the Gauss quadrature order-$N$ formula can also be employed
\beq
v_j(x)=\fr{\pi A_j}{N}\sum_{i=1}^N\fr{p(bx_i)}{|x/b-x_i|^\Ga}\sin\fr{(2i-1)\pi}{2N}, \quad x_i=\cos\fr{(2i-1)\pi}{2N},\quad |x|>b.
\label{5.10}
\eeq
The numerical tests show that in the case of Hertzian contact, when the pressure vanishes
at the endpoints, this approximation is in good agreement with the exact formulas  (A.5) and (A.6).

Consider the particular case $f(x)=Q_0x^2$. Owing to the fact that $g_n^{(1)}=0$ for all
$n$ except for $g_0^{(1)}$ and $g_2^{(1)}$  and employing formulas (\ref{4.14}) for these nonzero terms we
specify the formulas for the  parameter $\Gd$ and find explicitly the half-length of the contact zone
\beq
\Gd=\fr{Q_0b^2}{\Ga+2}
+\fr{\pi AP}{\GG_0 b\cos\fr{\pi\Ga}{2}},\quad
b=\left(\fr{\GG(2+\fr{\Ga}{2})\GG(\fr{1-\Ga}{2})(\Gt_1+\Gt_2)P}{\sqrt{\pi} Q_0}\right)^{\fr{1}{\Ga+2}}.
\label{5.11}
\eeq
For this parabolic case we also compute the pressure distribution and the normal displacement.
From (\ref{5.2}) we have
\beq
p(x)=\left(1-\fr{x^2}{b^2}\right)^{(\Ga-1)/2}
\left[\fr{P}{\GG_0 b}
+\fr{2b^2Q_0\cos\fr{\pi\Ga}{2}}{\pi A\Ga(\Ga+1)}\left(\fr{1}{\Ga+2}-\fr{x^2}{b^2}\right)\
\right].
\label{5.12}
\eeq
On substituting the expression (\ref{5.11}) for $b$ into formula (\ref{5.12}) we arrive at 
\beq
p(x)=\fr{P\GG(\fr{\Ga}{2}+2)}{\sqrt{\pi}b\GG(\fr{\Ga+3}{2})}\left(1-\fr{x^2}{b^2}\right)^{(\Ga+1)/2}.
\label{5.13}
\eeq
In the particular case, when one of the bodies is a rigid punch, this formula coincides with the corresponding expression of the pressure distribution derived by 
Giannakopoulos and  Pallot  (2000).  
 When $\Ga\to 0$, the contact zone half-length $b$
and the pressure distribution $p(x)$ tend to  $b_0$ and $p_0(x)$  which represent the contact zone half-length
and the contact pressure, respectively, in the case when both bodies are isotropic elastic bodies and whose contact is governed by equation (\ref{2.9})
$$
\lim_{\Ga\to 0} b=b_0, \quad b_0=\sqrt{\fr{(\Gt^\circ_1+\Gt_2^\circ)P}{Q_0}},
$$
\beq
\lim_{\Ga\to 0} p(x)=p_0(x), \quad p_0(x)=\fr{2P}{\pi b_0^2}\sqrt{b_0^2-x^2}.
\label{5.13'}
\eeq
where $\Gt_j^\circ$ are given by (\ref{2.8}).  This expression coincides with the contact pressure associated  with the elastic isotropic case of the Hertzian model and 
 obtained by solving equation (\ref{2.9}) (Shtayerman, 1949, Chapter II, (23)).

We determine the normal displacements of surface points $u_y(x,0)=-v_2(x)$ of the lower body $B_2$ outside the contact zone 
when $B_2$ is a half-plane that is when $f_2(x)=0$, while the upper body $B_1$ has the profile $y=Q_0x^2$. As before, the Young moduli of
both bodies have the 
same exponent and may have different factors $e_1$ and $e_2$. Since the only nonzero coefficients are $\GF_0$ and $\GF_2$,
we transform formula (\ref{5.8}) to the form
\beq
v_2(bt)=\Gt_2 b^{1-\Ga}[\GF_0 \tilde I_0(t;\Ga)+\GF_2 \tilde I_2(t;\Ga)],
\label{5.14}
\eeq
where
\beq
\GF_0=\fr{P}{b\GG_0}, \quad \GF_2=-\fr{4b^{\Ga+1}Q_0\cos\fr{\pi\Ga}{2}}{\pi\Ga(\Ga+1)(\Ga+2)(\Gt_1+\Gt_2)}, \quad 
\tilde I_n(t;\Ga)=\fr{1}{\Ga} I_n(t;\Ga). 
\label{5.15}
\eeq
From (A.5),
$$
\tilde I_0(t;\Ga)=\fr{\pi}{\cos\fr{\pi\Ga}{2}}\left[\fr{1}{\Ga}-\fr{\GG(\fr{\Ga+1}{2})(-\fr{t+1}{2})^{(1-\Ga)/2}}{\GG(\Ga+1)\GG(\fr{3-\Ga}{2})}
F\left(\fr{1-\Ga}{2},\fr{1+\Ga}{2},\fr{3-\Ga}{2};\fr{t+1}{2}\right)\right],
$$
$$
\tilde I_2(t;\Ga)=\fr{\pi\Ga(\Ga+1)}{2\cos\fr{\pi\Ga}{2}}\left[\fr{\Ga+1}{2}F\left(-2,\Ga+2,\fr{\Ga+1}{2};\fr{t+1}{2}\right)
\right.
$$
\beq
\left.
-\fr{\GG(\fr{\Ga+1}{2})(-\fr{t+1}{2})^{(1-\Ga)/2}}{\GG(\Ga+1)\GG(\fr{3-\Ga}{2})}
F\left(-\fr{3+\Ga}{2},\fr{5+\Ga}{2},\fr{3-\Ga}{2};\fr{t+1}{2}\right)\right],\quad -3<t<-1.
\label{5.16}
\eeq
For $t<-3$ the integrals $\tilde I_n(t;\Ga)=\Ga^{-1} I_n(t;\Ga)$  ($n=0,2$) are obtained from formula (A.6). 
It is easy to see from formulas (\ref{5.14}) to (\ref{5.16}) that the displacements $v_j(t;\Ga)$ become 
infinite when $\Ga\to 0$, and the limit transition  $\Ga\to 0$ for the displacements outside the contact zone is impossible.

\begin{figure}[t]
\centerline{
\scalebox{0.6}{\includegraphics{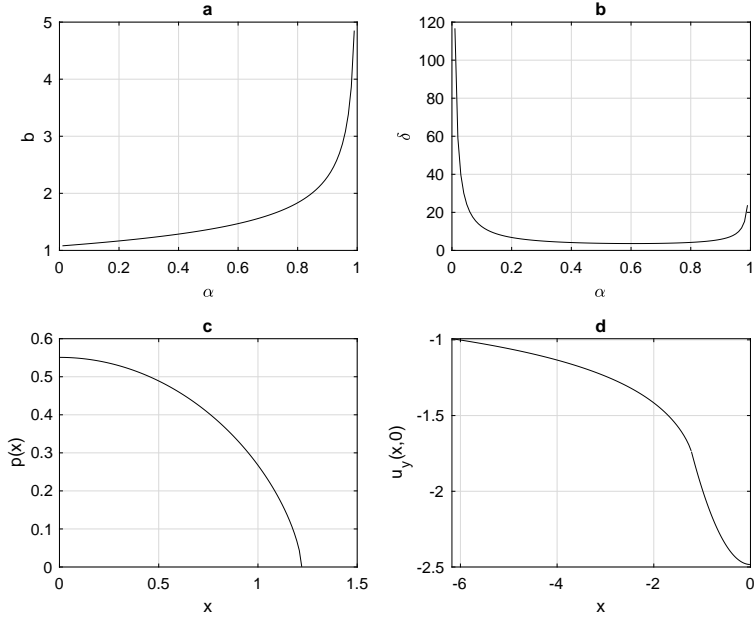}}}
\caption{
 The case $\Ga_1=\Ga_2=\Ga$ when $P=1$, $f(x)=x^2$, $e_1=e_2=1$, $\nu_1=\nu_2=0.3$. (a): the contact zone 
 half-length $b$ versus $\Ga\in(0,1)$, 
 (b): the parameter $\Gd$ versus $\Ga$, (c): pressure $p(x)$ for $\Ga=0.3$, (d): the displacement $u_y(x,0)=-v_2(x)$ for $x<-b$
 and $u_y(x,0)=(x^2-\Gd)/2$ for $x\in(-b,0)$ when
 $\Ga=0.3$.}
\label{fig8}
\end{figure}

Figure 8  shows the results of computations in the case when the Young moduli of the contacting bodies are the same,
$E_j(z)=e_j z^\Ga$ and $e_1=e_2$. We choose
$P=1$, $f(x)=x^2$, $e_1=e_2=1$, and $\nu_1=\nu_2=0.3$.
Figures 8 (a) and (b) demonstrate the variation of the half-length $b$ and the parameter $\Gd$ 
with the exponent $\Ga\in(0,1)$. 
 As $\Ga$ grows the contact zone becomes larger. As in the
case $\Ga_1\ne \Ga_2$, as $\Ga\to 0$, the parameter $\Gd\to\infty$. It also grows as $\Ga\to 1$. 
In Figures 8 (c)  and (d), we present sample curves for the contact pressure for $x\in[0,b]$ and the normal displacement 
$u_y(x,0)=-v_2(x)$ of the surface of the lower body (a half-plane) when $f_1(x) =x^2$, $f_2(x)=0$ for $x<-b$
and  $u_y(x,0)=(x^2-\Gd)/2$ for $x\in(-b,0)$. In both Figures 8 (c) and 8 (d),   
$\Ga=0.3$ (in this case $b= 1.22072$).

\setcounter{equation}{0}

\section{Surface energy model}

In this section following the JKR model (Johnson et al, 1971; Johnson, 1985) we aim to take into account the effect of adhesive forces (Figure 1 (b)) and study their impact
on the contact zone size, the contact pressure and the normal displacement. In the two-dimensional
case, the loss of surface energy is given by  $U_s=-2\Gg_s b$, where $\Gg_s$ is the work of adhesion (a half-density
of the surface energy). The elastic strain energy is expressed through the normal displacement
$v(x)=\Gd-f(x)$, and the contact pressure $p(x)$ as
\beq
U_e=\fr12\int_{-b}^b p(x) v(x)dx,
\label{6.1}
\eeq
and the total energy defined by $U_{total}=U_e-2\Gg_s b$ is a function of the contact zone half-length $b$. 
In contrary to the Hertzian model, the JKR model admits singularities of the contact pressure at the endpoints.
Also, the parameter $b$ is defined not from the condition that quenches  the pressure singularities but from
the condition of minimum of the total energy that is
\beq
\fr{d U_e}{db} -2\Gg_s=0.
\label{6.2}
\eeq

\subsection{Case $\Ga_1=\Ga_2=\Ga$}

\begin{figure}[t]
\centerline{
\scalebox{0.6}{\includegraphics{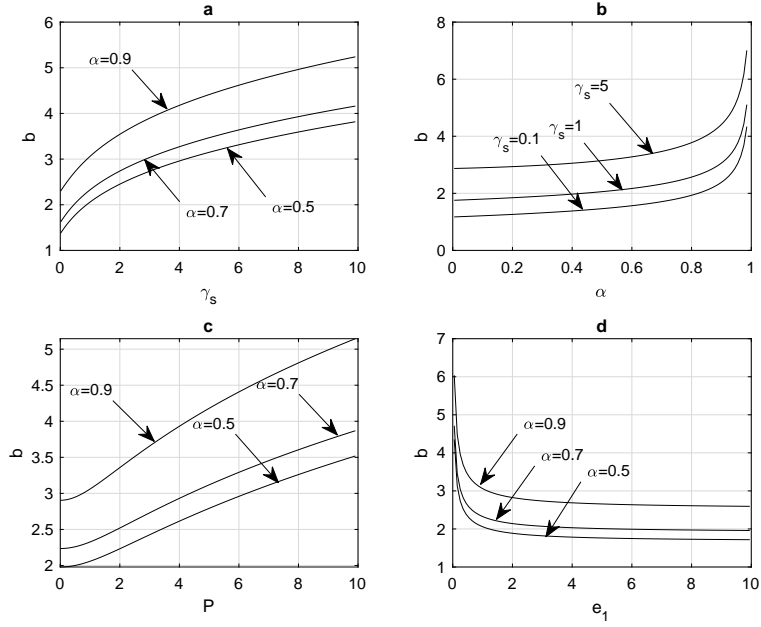}}}
\caption{
JKR-model: the contact zone half-length $b$ in the case $\Ga_1=\Ga_2=\Ga$, $f(x)=x^2$, $e_2=1$, $\nu_1=\nu_2=0.3$.
 (a): $b$ versus the surface energy density $\Gg_s$, 
 (b):  $b$ versus the exponent   $\Ga$, (c): $b$ versus the normal force $P$ when $\Gg_s=1$, 
 (d):   $b$ versus the parameter $e_1$  when $\Gg_s=1$.}
\label{fig9}
\end{figure}

\begin{figure}[t]
\centerline{
\scalebox{0.6}{\includegraphics{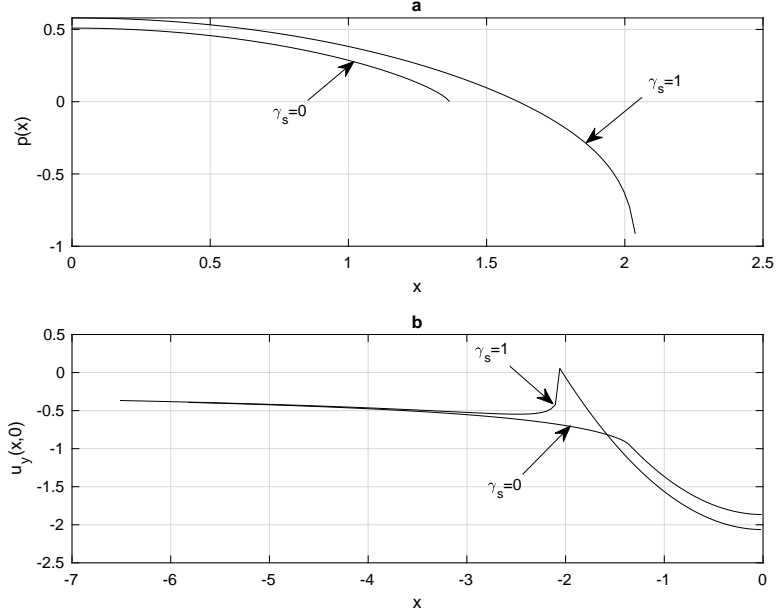}}}
\caption{(a): the contact pressure $p(x)$ and the normal displacement $u_y(x,0)=-v_2(x)$ for $x<-b$
 and $u_y(x,0)=(x^2-\Gd)/2$ for $x\in(-b,0)$ for the Hertzian ($\Gg_s=0$) and JKR ($\Gg_s=1$) models
when  $\Ga=0.5$, $P=1$, $f(x)=x^2$, $e_1=e_2=1$, $\nu_1=\nu_2=0.3$.}
\label{fig10}
\end{figure}

We consider  parabolic profiles  of the contacting bodies, $f(x)=Q_0x^2$.
In this case,  the pressure $p(x)$ is given by (\ref{5.12}) and has order $(\Ga-1)/2$ power singularities at the endpoints, while the 
parameter $b$ is free.  Since the resultant force $P$  has to be balanced by the contact pressure $p(x)$,
we satisfy the condition (\ref{2.2}) and   define $\Gd$ by formula (\ref{5.11}). 
To evaluate the integral (\ref{6.1}) we write the pressure
and displacement in the following equivalent form:
$$
p(bt)=(1-t^2)^{(\Ga-1)/2}[\GF_0 C_0^{\Ga/2}(t)+\GF_2 C_2^{\Ga/2}(t)], \quad 
$$
\beq
v(bt)=\Gd C_0^{\Ga/2}(t)-\fr{Q_0 b^2}{\Ga(\Ga+2)}[\Ga C_0^{\Ga/2}(t)+2 C_2^{\Ga/2}(t)],
\label{6.3}
\eeq
where $\GF_0$ and $\GF_2$ are determined in (\ref{5.15}). Using the orthogonality property (\ref{3.10}) of the Gegenbauer polynomials we obtain
\beq
U_e=\fr{P^2(\Gt_1+\Gt_2)\GG(\fr{\Ga}{2}+1)\GG(\fr{1-\Ga}{2})b^{-\Ga}}{2\Ga\sqrt{\pi}}+
\fr{\sqrt{\pi} Q_0^2 b^{\Ga+4}}{2(\Gt_1+\Gt_2)(\Ga+2)\GG(\fr{1-\Ga}{2})\GG(\fr{\Ga}{2}+3)}.
\label{6.4}
\eeq
The derivative of the   elastic strain energy  in (\ref{6.2}) can be evaluated exactly, and we arrive at the following 
transcendental equation with respect to the parameter $b$:
\beq
\fr{\sqrt{\pi} Q_0^2 b^{2\Ga+4}}{(\Gt_1+\Gt_2)(\Ga+2)\GG(\fr{1-\Ga}{2})\GG(\fr{\Ga}{2}+2)}-2\Gg_s b^{\Ga+1}-
\fr{P^2(\Gt_1+\Gt_2)\GG(\fr{\Ga}{2}+1)\GG(\fr{1-\Ga}{2})}{2\sqrt{\pi}}=0.
\label{6.5}
\eeq
Passing to the limit $\Ga\to 0$ and keeping $\Gg_s\ge 0$ we reduce the transcendental equation to the quartic equation
\beq
\fr{Q_0^2 b^{4}}{2(\Gt^\circ_1+\Gt^\circ_2)}-2\Gg_s b-
\fr{P^2(\Gt^\circ_1+\Gt^\circ_2)}{2}=0,
\label{6.6}
\eeq
and when, in addition, $\Gg_s\to 0$, we obtain the classical formula (\ref{5.14}) for the value of $b$ in the case of Hertzian contact
of two elastic isotropic bodies.

Passing to the limit $\Gg_s\to 0$ in equation (\ref{6.5})  and keeping $\Ga\in(0,1)$ we arrive at the equation with respect to $b$
that admits an exact solution; it coincides with the value of $b$ in the Hertzian model given by (\ref{5.11}).

If we assume that $\nu_1=\nu_2=0.5$, by passing to the limit $\Ga\to 1$ 
we derive from (\ref{6.5}) the following cubic equation with respect to $b^2$ for two Gibson solids (Gibson, 1967):
\beq
\fr{8}{27}Q_0^2 b^6-2\Gg_s\left(\fr{1}{e_1}+\fr{1}{e_2}\right)b^2-\fr{3P^2}{8}\left(\fr{1}{e_1}+\fr{1}{e_2}\right)^2=0.
\label{6.7}
\eeq

Notice that the transcendental equation (\ref{6.5})  is different from the corresponding equations obtained by  
Giannakopoulos and Pallot (2000) and  Chen et al (2009a). These authors split the solution into two parts, the first one gives the solution
for a parabolic punch and the second one corresponds to the model of a flat punch. 
The discrepancies between the transcendental equations obtained by these authors and equation (\ref{6.5})
are caused by their disregard for the fact that
 the displacement $\Gd$  is a function of  the contact zone half-length $b$. 
This explains why the limit transition $\Ga\to 0$ is
impossible in the solutions obtained by these authors. 

A number of tests have been conducted to ascertain  the impact of the surface energy density $2\Gg_s$ and the exponent $\Ga$
on the contact zone size $2b$, the pressure distribution, and the normal displacement.  The curves in Figure 9 (a) exhibit
an increase of the contact zone size with the parameter $\Gg_s$. It is seen from Figure 9 (b) that the half-length $b$ rapidly
increases when the exponent $\Ga$ approaches  1.  The $P-b$ curves  in Figure 9 (c) demonstrate the rate of growth
of the half-length $b$  when the total force $P$ grows. From  Figure 9 (d) it is seen that  when the Young moduli are
$E_j=e_j |y|^\Ga$, $e_2=1$ and the factor $e_1$ grows, the parameter $b$ first rapidly decreases and then 
its rate of decrease is insignificant. 

Contact pressure curves computed according to the Hertz  and JKR models are portrayed in Figure 10 (a). Since the pressure
$p(x)$ is an even function,
the curves demonstrate that the pressure vanishes at the endpoints $x=\pm b$ in the former model. In the JKR model,
the contact stress is compressive for $-b_*<x<b_*$ and tensile at the edge zones $(-b,-b_*)$ and $(b_*,b)$. 
The numerical tests show that a growth of the surface energy density $2\Gg_s$ shrinks  the central zone, where the stress is compressive,
and enlarges the zone, where the stress is tensile.

As in the Hertzian contact model, we analyze the  normal displacements $u_y(x,0)=-v_2(x)$ for the JKR model outside of the contact zone 
when the upper body has a parabolic profile, $f_1(x)=Q_0 x^2$, and the lower body is a half-plane, $f_2(x)=0$.
The displacement $v_2(x)$ is given by the same formula (\ref{5.14}).  However, since the pressure $p(x)$ does not vanish
at the endpoints, has power singularities  and is described by formula  (\ref{5.12}), the derivative of the displacement (\ref{5.14})
tends to infinity as $x\to\pm b$. A sample curve of the displacement $u_y(x,0)=-v_2(x)$ for $x<-b$ 
 and $u_y(x,0)=(x^2-\Gd)/2$ for $x\in(-b,0)$  ($Q_0=1$)  is shown in Figure 10 (b). It is seen that in the case of Hertzian
 contact ($\Gg_s=0)$ the contact surface is smooth near the contact zone endpoints, while in the case of the JKR model, due
 to the adhesion forces a part of the surface of the flat body $B_2$ is attracted to the interface, and the tangent lines
 to the surfaces of the contacting bodies at the endpoints have different slopes.

\subsection{Case $\Ga_1>\Ga_2$}

The direct method for computing the  elastic strain energy described in the previous section can be generalized to the case
when the contacting bodies have different exponents and as before, $\Ga_1>\Ga_2$. 
On expanding the normal displacement $v(x)=\Gd-f(x)$ $(-b<x<b)$ through the Gegenbauer polynomials of even order 
we have
\beq
v(bt)=\Gd-\sum_{k=0}^\infty a_{2k} C_{2k}^{\Ga_1/2}(t), \quad -1<t<1,
\label{6.8}
\eeq
substituting it together with the contact pressure 
\beq
p(bt)=(1-t^2)^{(\Ga_1-1)/2}\sum_{n=0}^\infty[\GF_{2n}^{(1)}+\Gd\GF_{2n}^{(2)}]C_{2n}^{\Ga_1/2}(t), \quad -1<t<1, 
\label{6.9}
\eeq
into formula (\ref{6.1}) and using the orthogonality of the polynomials we derive the series representation
of the elastic strain energy
\beq
U_e=\fr{b\Gd}{2}
(\GF_0^{(1)}+\Gd\GF_0^{(2)})\GG_0-\fr{b}{2}\sum_{n=0}^\infty
(\GF_{2n}^{(1)}+\Gd\GF_{2n}^{(2)})a_{2n} h_{2n}(\Ga_1).
\label{6.10}
\eeq
As before, we simplify the formula for the parabolic case, $f(x)=Q_0x^2$. In this case, $a_n=0$ unless $n=0$ or $n=2$,
\beq
a_0=\fr{b^2 Q_0}{\Ga_1+2}, \quad a_2=\fr{2 b^2 Q_0}{\Ga_1(\Ga_1+2)},
\label{6.11}
\eeq
and ultimately we have
\beq
U_e=\fr{b\GG_0}{2}(\GF_0^{(1)}+\Gd\GF_0^{(2)})\left(\Gd-\fr{b^2Q_0}{\Ga_1+2}\right)
-\fr{2b^3Q_0\sqrt{\pi}\GG(\fr{\Ga_1+3}{2})}
{(\Ga_1+2)(\Ga_1+4)\GG(\fr{\Ga_1}{2})}(\GF_{2}^{(1)}+\Gd\GF_{2}^{(2)}).
\label{6.12}
\eeq
The minimum of the total energy attaines if the contact zone half-length $b$ solves the  transcendental equation
\beq
\fr{dU_e}{db}-2\Gg_s=0.
\label{6.13}
\eeq
Explicit differentiation is impossible for the coefficients $\GF_0^{(j)}$ and  $\GF_2^{(j)}$
being a part of the solution to the infinite system (\ref{3.36}), and there is no way to explicitly separate $b$ 
from the unknowns of the 
infinite system. Approximately, equation (\ref{6.13}) can be written as
\beq
 \fr{U_e(b+\Gve)-U_e(b)}{\Gve}-2\Gg_s\approx 0, 
\label{6.14}
\eeq
where $\Gve$ is a small and positive.

\begin{figure}[t]
\centerline{
\scalebox{0.6}{\includegraphics{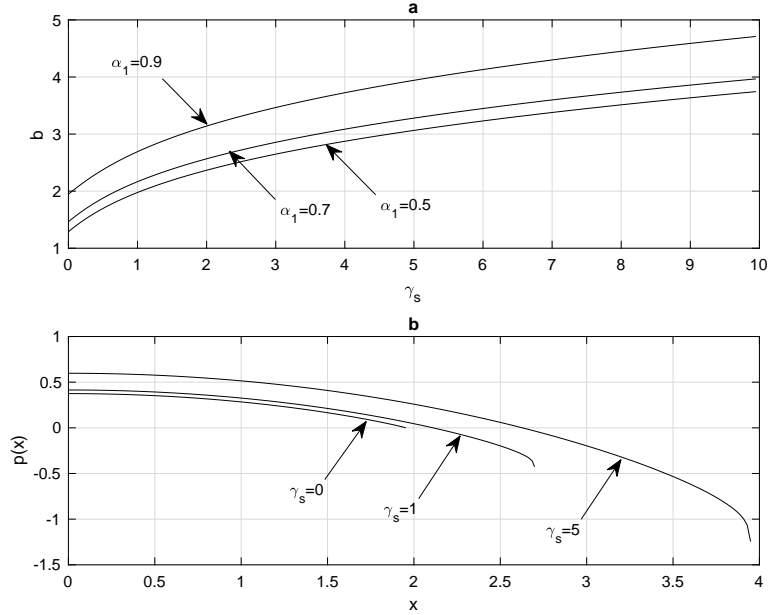}}}
\caption{JKR model for the case of different exponents: $E_1=e_1 y^{\Ga_1}$ and $E_2=e_2 (-y)^{\Ga_2}$.
when $e_1=e_2=1$, $P=1$, $Q_0=1$, $\nu_1=\nu_2=0.3$.
(a): the half-length $b$ versus the  half-density $\Gg_s$ of the surface energy for  
$\Ga_1=0.5, 0.7, 0.9$ and $\Ga_2=\Ga_1/2$.
(b): the contact pressure $p(x)$ for $x\in(0,b)$ for $\Gg_s=0, 1, 5$ when $\Ga_1=0.9$ and $\Ga_2=0.5$.}
\label{fig11}
\end{figure} 
The variation of the half-length of the contact zone $b$ with the half-density $\Gg_s$ of surface energy 
for three values of the exponent $\Ga$ is portrayed in Figure 11 (a). It has been calculated by the method of orthogonal polynomials
presented in Section 3. The difference between the scheme for the Hertzian and JKR models is only in the way 
how the parameter $b$ is fixed. In the Hertzian model, it solves the transcendental equation (\ref{4.8}) that guarantees
that the pressure vanishes at the endpoints, while in the JKR model, it is defined from the approximate
equation (\ref{6.14}), the condition of minimum of the total energy. 
For computations, $\Gve$ is accepted to be $10^{-4}$, and the differences between the results for $\Gve=10^{-3}, 10^{-4}, 10^{-5}$
are not significant. For example, for $\Ga_1=0.5$, $\Ga_2=0.25$, $e_1=e_2=1$, $P=1$, $Q_0=1$, $\nu_1=\nu_2=0.3$,
and $\Gg_s=1$ we have $b=1.97621$ if  $\Gve=10^{-3}$, $b=1.97666$ if $\Gve=10^{-4}$, and $b=1.97670$ if  $\Gve=10^{-5}$.
It turns out that as $\Gg_s\to 0$, the parameter
$b$ associated with the JKR model tends to the one for the Hertzian model. It is not seen how this result can be proved
analytically. However, all our numerical tests confirm this conclusion.

The pressure distribution $p(x)$ is shown in Figure 11 (b) for three values of the parameter $\Gg_s$ when $\Ga_1=0.9$
and $\Ga_2=0.5$. As $\Gg_s\to 0$,  the contact pressure vanishes at the endpoints and coincides with the pressure
found from the Hertzian model. When $\Gg_s>0$, similarly to the case $\Ga_1=\Ga_2$, the contact zone enlarges and the contact stress becomes tensile at two edge
zones $(-b,-b_*)$ and $(b_*,b)$ and tends to $\infty$ as $x\to \pm b$. 

\vspace{.1in}

{\large {\bf Conclusions}}

\vspace{.1in}

We analyzed two plane problems, the Hertzian and JKR models,  of frictionless contact of two inhomogeneous elastic bodies with distinct
moduli of elasticity   $E_1(y)=e_1 y^{\Ga_1}$ and $E_2(y)=e_2 (-y)^{\Ga_2}$ with $0<\Ga_2\le\Ga_1<1$.
On employing the Rostovtsev representation of the normal displacement in the contact zone through the pressure
distribution we showed that the model is governed by an integral equation with two different power kernels.
For its solution, a novel method of Gegenbauer orthogonal polynomials was proposed. We  reduced the integral
equation to an infinite system of linear algebraic equations whose coefficients, after some transformations, become
integral free. It was demonstrated that when $\Ga_2\to\Ga_1$, the infinite system is decoupled, and its exact solution
coincides with the one obtained by direct solution of the integral equation with one power kernel.

We found a rigid body displacement  $\Gd$ (the total displacement of distant points of the bodies) from the equilibrium condition that balances the normal
total force and the contact pressure.   The length of the contact zone is determined from a
transcendental equation
that guarantees that the pressure vanishes at the endpoints in the Hertzian model and the total energy attains its minimum
in the JKR model. The pressure distribution is found in a series form,  and the coefficients of the 
expansion are determined from an infinite system of the second kind solved by the reduction method. 
The numerical tests implemented
revealed    rapid convergence of the method for all admissible values of the model parameters.
By employing the method of Mellin's convolution integrals and the theory of residues, we computed the normal displacements of surface points outside the contact zone. In the Hertzian model, the profile of the contacting surfaces at the endpoints is smooth, while
in the JKR model, the derivative of the normal displacement is infinite, and  a part of the contacting surfaces 
is attracted by adhesion forces to the interface. In contrary to the Hertzian model, the pressure distribution
does not vanish at the endpoints, it tends to $-\infty$, and there are two edge zones in the contact area where the contact stress is tensile.  

Our numerical results showed that the parameter $\Gd$, the contact zone length, the contact pressure, and the elastic displacement
significantly depend on variation of the bigger parameter $\Ga_1$ and only slightly vary with the second, smaller, parameter $\Ga_2$.
When the exponent $\Ga_1$ is growing, the contact zone is also growing. 
In the case when the two exponents $\Ga_1$ and $\Ga_2$ are the same, $\Ga_1=\Ga_2=\Ga$, we obtained the contact zone
length, the parameter $\Gd$, the contact pressure and the normal displacements of the surface points exactly.
By passing to the limit $\Ga\to 0$, we showed that the result coincides with the classical solution of the problem of Hertzian
contact of two isotropic elastic bodies. 

For  the JKR model, we found out that when the half-density of surface energy $\Gg_s\to 0$,
the contact zone length, pressure and normal displacement tend to the corresponding quantities associated with the Hertzian model.
In both cases, $\Ga_1>\Ga_2$ and $\Ga_1=\Ga_2$, the transcendental equation for the contact zone length admits 
passing to the limit $\Ga_j\to 0$. This is possible not only for the contact zone length but also for the contact pressure. However,
the normal displacement derived for both Hertzian and JKR  models  when $\Ga_1\ge \Ga_2>0$ become infinite
when $\Ga_1\to 0$. This is due to the presence of $\Ga_1^{-1}$ in the formula for the displacement $\Gd$ of distant points   
of the contacting bodies.

The method we presented admits generalizations and modifications in different directions including the Hertzian and JKR axisymmetric contact models
of two power-law graded bodies. 

%\section*
%{\large {\it  Acknowledgments.}} This work was partly funded by NSF through grant DMS0707724.

\vspace{.1in}

{\large {\bf Appendix A. Evaluation of the integral $I_n(t;\Ga_j)$}

\vspace{.1in}

We wish to evaluate the integral
$$
 I_{n}(t;\Ga_j)=\int_{-1}^1\fr{(1-\Gj^2)^{(\Ga_1-1)/2} C_n^{\Ga_1/2}(\Gj)d\Gj}{|\Gj-t|^{\Ga_j}}, \quad 0<\Ga_j<1, \quad t<-1.
$$
 First, we  transform the integral  to the form
$$
I_{n}(t;\Ga_j)=2^{\Ga_1-\Ga_j}\int_0^1\fr{[\Gn(1-\Gn)]^{(\Ga_1-1)/2} C_n^{\Ga_1/2}(2\Gn-1)d\Gn}
{(\Gn+\Gz)^{\Ga_j}},
$$
where $\Gj=2\Gn-1$, $t=-2\Gz-1$, and $\Gz>0$.
Next, we represent the integral $I_{n}(t;\Ga_j)$ as a Mellin convolution integral
$$
I_{n}(t;\Ga_j)=2^{\Ga_1-\Ga_j}\int_0^\infty h_1(\Gn)h_{2}\left(\fr{\Gz}{\Gn}\right)
\fr{d\Gn}{\Gn},
$$
where
$$
h_1(\Gn)=\left\{
\begin{array}{cc}
\Gn^{(1+\Ga_1)/2-\Ga_j}(1-\Gn)^{(\Ga_1-1)/2}C_n^{\Ga_1/2}(2\Gn-1), & 0<\Gn<1\\
0, & \Gn>1\\
\end{array}
\right.,
\quad h_{2}(\Gz)=\fr{1}{(1+\Gz)^{\Ga_j}}.
$$
We aim further to apply the Mellin convolution theorem (Titchmarsh,  1948)
$$
I_n(t;\Ga_j)=\fr{2^{\Ga_1-\Ga_j}}{2\pi i}\int_{\Gs_j-i\infty}^{\Gs_j+i\infty}H_1(s)H_{2}(s)\Gz^{-s}ds,
\eqno(A.3)
$$
where $H_1(s)$ and $H_{2}(s)$ are the Mellin transforms of the functions $h_1(\Gn)$ and $h_{2}(\Gn)$, respectively.
These transforms are obtained by exploiting the following integrals (Gradshteyn and Ryzhik, 1994,  formulas 7.311(3) and 3.194(3)):
$$
H_1(s)=\fr{\GG(\fr{\Ga_1+1}{2})(\Ga_1)_n}{(-1)^nn!}
\fr{\GG(s+\fr{\Ga_1+1}{2}-\Ga_j)\GG(-s+\Ga_j+n)}{\GG(s+\Ga_1-\Ga_j+n+1)\GG(-s+\Ga_j)},\quad \R s>\Ga_j-\fr{\Ga_1+1}{2},
$$
$$
H_{2}(s)=\fr{\GG(s)\GG(\Ga_j-s)}{\GG(\Ga_j)},\quad 0<\R s<\Ga_j,
\eqno(A.4)
$$
and since $\Ga_j\in (0,1)$ and $\Ga_1>\Ga_2$,  we have $\Gs_j \in(0,\Ga_j)$. The final step of the procedure is substituting formulas (A.4) into
(A.3) and applying the theory of residues. In the case $0<\Gz<1$ ($-3<t<-1$) this implies 
$$
I_n(t;\Ga_j)=\fr{(\Ga_1)_n \GG(\fr{\Ga_1+1}{2})}{(-1)^n 2^{\Ga_j-\Ga_1}n!}
\left[
\fr{(\Ga_j)_n\GG(\fr{\Ga_1+1}{2}-\Ga_j)}{\GG(\Ga_1-\Ga_j+n+1)}
 F(\Ga_j-\Ga_1-n, \Ga_j+n, \Ga_j+\fr{1-\Ga_1}{2}; -\Gz)
\right.
$$
$$
\left.
+\fr{\GG(\Ga_j-\fr{\Ga_1+1}{2})}{\GG(\Ga_j)}\Gz^{(\Ga_1+1)/2-\Ga_j}
F\left(\fr{1+\Ga_1}{2}+n,\fr{1-\Ga_1}{2}-n,\fr{3+\Ga_1}{2}-\Ga_j; -\Gz\right)
\right], \quad 0<\Gz<1,
\eqno(A.5)
$$
where $F$ is the hypergeometric function. If $\Gz>1$, that is if $t<-3$, then we have 
$$
I_{n}(t;\Ga_j)=\fr{(-1)^n\sqrt{\pi}(\Ga_1)_n(\Ga_j)_n\GG(\fr{\Ga_1+1}{2})}{2^{\Ga_j+2n}n!\GG(\fr{\Ga_1}{2}+n+1)\Gz^{\Ga_j+n}}
 F\left(\Ga_j+n,\fr{\Ga_1+1}{2}+n, \Ga_1+2n+1; -\fr{1}{\Gz}\right).
\eqno(A.6)
$$
In a neighborhood of the point $\Gz=1$ ($t=-3$), for computational purposes, it is numerically efficient to employ the formula
9.131(1) (Gradshteyn and Ryzhik, 1994) that is 
$$
F(\Ga,\Gb,\Gg;z)=(1-z)^{-\Gb}F\left(\Gb,\Gg-\Ga,\Gg;\fr{z}{z-1}
\right).
$$

\vspace{.2in}

{\large {\bf References}}

\vspace{.1in}

Antipov, Y.A.,  Mkhitaryan, S.M., 2021.  Integral and integro-differential equations with an exponential kernel and applications. 
Quart. J. Mech. Appl. Math.  74,  297-322. 

Bateman, H. 1954. Tables of Integral Transforms, vol. 2. Bateman Manuscript Project, McGraw-Hill, New York.

Chen, S., Yan, C., Soh, A., 2009a. Adhesive behavior of two-dimensional power-law
graded materials. Int. J. Solids Struct. 46, 3398–3404.

Chen, S., Yan, C., Zhang, P., Gao, H., 2009b. Mechanics of adhesive contact on a power-law graded elastic half-space. J. Mech. Phys. Solids 57, 1437–1448.

 Gakhov, F.D., 1966. Boundary Value Problems. Pergamon Press, Oxford.

Giannakopoulos, A.E., Suresh, S., 1997a. Indentation of solids with gradients in elastic properties: part I. Point force. Int. J. Solids Struct. 19, 2357–2392.

Giannakopoulos, A.E., Suresh, S., 1997b. Indentation of solids with gradients in elastic properties: part II. Axisymmetric indentors. Int. J. Solids Struct. 34,
2393–2428.

Giannakopoulos, A.E., Pallot, P., 2000. Two-dimensional contact analysis of elastic
graded materials. J. Mech. Phys. Solids 48, 1597–1631.

Gibson, R.E., 1967. Some results concerning displacements and stresses in a non-homogeneous elastic half-space. Geotechnique 17, 58–67.

Gradshtein, I.S., Ryzhik, I.M., 1994. Table of integrals, series, and products. Academic Press, New York.

Guo, X., Jin, F.,  Gao, H.,  2011. Mechanics of non-slipping adhesive contact on a power-law graded elastic
half-space. Int.  J. Solids Struct. 48, 2565–2575.

Gutleb, T.S.,  Carrillo, J.A., Olver, S. 2021. Computing equilibrium measures with power law kernels. arXiv:2011.00045.

Jin, F.,  Tang  Q., Guo, X., Gao, H.,  2021.
A generalized Maugis-Dugdale solution for adhesion of power-law
graded elastic materials. J. Mech. Phys. Solids 154, 104509.

Johnson, K.L.,  1985. Contact Mechanics. Cambridge University Press, Cambridge.

Johnson, K.L.,   Kendall, K.,  Roberts, A.D., 1971. Surface energy and the contact of elastic solids. Proc Roy. Soc. A 324,
301-313.

Klein, G.K., 1956.   Allowing for inhomogeneity, discontinuity of the deformations and other mechanical properties of the 
soil  in the design of structures on a continuous foundation. Sb. Trudov Mosk. (Moscow)  Inzh.-Str. Inst. 14, 168-180.

Korenev, B. G., 1957. A die resting on an elastic half-space, the modulus of elasticity of which is an
exponential function of depth. Dokl. Akad. Nauk SSSR, 112, 823–826.

Korenev, B.G., 1960. Some Problems of the Theory of Elasticity and Heat Conduction Solvable by Bessel Functions.
Fizmatgiz, Moscow.

Lekhnitskii, S.G., 1962. Radial distribution of stresses in a wedge and in a half-plane with variable modulus of elasticity.
J. Appl. Math. Mech. (PMM) 26, 199-206.

Maugis, D., 1992. Adhesion of spheres: The JKR-DMT transition using a Dugdale model. J. Colloid Interface Sci. 150, 243–269.

Mkhitaryan, S.M.,  2015. An eigenvalue relation in spheroidal wave functions related to
potential theory and its applications to contact problems,  J. Appl. Math. Mech. (PMM) 79, 304–313. 

Mossakovskii, V.I., 1958. Pressure of a circular die [punch] on an elastic half-space, whose modulus of elasticity is an exponential [power] function of depth.  J. Appl. Math. Mech. 
(PMM) 22,  168-171. 

Popov, G.Ia., 1961.  On a method of solution of the axisymmetric contact problem of the theory of elasticity.
J. Appl. Math. Mech. (PMM)  25, 105–118.

Popov, G.Ya., 1967. On an approximate method of solution of a contact problem of an annular punch. 
Izv. AN Arm SSR, Mekhanika 20, 2.

Popov, G.Ya., 1971. Contact problem of elasticity when there is a circular contact region and the surface structure of the contacting
bodies is taking into account. Izv. AN SSSR, Mekh. Tv. Tela, 3. 

Popov, G.Ia., 1973. Axisymmetric contact problem for an elastic inhomogeneous half-space in the presence of cohesion. J. Appl. Math. Mech. 
(PMM) 37,  1052-1059.

Popov, G.Ya., 1982. Concentration of elastic stresses near stamps, cuts, thin inclusions, and
reinforcements. Nauka, Moscow.

Rostovtsev, N.A., 1964. On the theory of elasticity of a nonhomogeneous medium. J. Appl. Math. Mech. (PMM) 28, 745-757.

Saleh, B.,  Jiang, J., Fathi, R., Al-hababi,  T.,  Xu, Q., 
Wang, L.,  Song, D.,  Ma, A.,  2020.
30 Years of functionally graded materials: An overview of manufacturing methods, applications and future challenges.
Composites Part B: Eng. 201, 108376.

Shtayerman, I. Ya., 1949. Contact Problem of the Theory of Elasticity. Gostekhizdat, Moscow. (Engl. Transl.: FTD-MT- 24-61-70
 by  Foreighn Technology Division, WP-AFB, Ohio.)

Titchmarsh, E.C. 1948. Introduction to the Theory of Fourier Integrals, Clarendon Press, Oxford.

Whipple, F.J.W., 1925. A group of generalized hypergeometric series: relations between 120 allied series of the type $F\left[\begin{array}{ccc}
a, & b, & c\\ & e, & f\\ \end{array}\right]$. Proc. London Math. Soc. 23, 104-114.

Willert, E., 2018. Dugdale-Maugis adhesive normal contact of axisymmetric power-law graded elastic bodies. Facta Univers. Ser. Mech. Eng. 16, 9–18.

\end{document}